\documentclass{amsproc}

\usepackage{enumerate}
\usepackage{mathtools}
\usepackage{tikz}
\usepackage{mathrsfs}
\usepackage{amsmath}
\usepackage{amssymb}
\usepackage{tikz-3dplot}
\usepackage{youngtab}

\theoremstyle{definition}
\newtheorem{theorem}{Theorem}[section]

\newtheorem{definition}[theorem]{Definition}
\newtheorem{example}[theorem]{Example}

\newtheorem{conjecture}[theorem]{Conjecture}
\newtheorem{question}[theorem]{Question}
\newtheorem{proposition}[theorem]{Proposition}
\newtheorem{corollary}[theorem]{Corollary}

\theoremstyle{remark}
\newtheorem{remark}[theorem]{Remark}

\numberwithin{equation}{section}

\def\PP{\mathbb{P}}
\def\RR{\mathbb{R}}

\def\ZZ{\mathbb{Z}}
\def\QQ{\mathbb{Q}}

\def\L{\mathbb{L}}

\DeclareMathOperator{\vol}{vol}

\begin{document}

\title{Volume polynomials}


\author{June Huh}
\address{Princeton University and Korea Institute for Advanced Study}
\email{huh@princeton.edu}
\thanks{The author is  partially supported by the Simons Investigator Grant.}


\subjclass[2020]{Primary 14-06. Secondary 52-06.}

\date{}

\begin{abstract}
Volume polynomials form a distinguished class of log-concave polynomials with remarkable analytic and combinatorial properties. 
I will  survey realization problems related to them, review fundamental inequalities they satisfy, and discuss applications to the combinatorics of algebraic matroids.
These notes are based on lectures given at the 2025 Summer Research Institute in Algebraic Geometry at Colorado State University.
\end{abstract}

\maketitle


\section{Realization problems for projection volumes and homology classes}

\subsection{}\label{SubsectionConvex}

Let $\pi_{ij}$ be the coordinate projection of $\mathbb{R}^4$ onto the plane orthogonal to the standard basis vectors $\mathbf{e}_i$ and $\mathbf{e}_j$.
For a convex body $A$ in $\mathbb{R}^4$, we consider its vector of projection areas $(p_{12}, p_{13}, p_{14}, p_{23}, p_{24}, p_{34})$, where
\[
p_{ij}=\big(\text{the area of the projection $\pi_{ij}(A)$}\big).
\]
Which tuples of six nonnegative real numbers can arise in this way?
This question is the simplest nontrivial instance of the various \emph{realization problems} for volume polynomials (Sections~\ref{SectionConvex} and ~\ref{SectionProjective}).
The following answer was given in  \cite[Theorem 1.4]{HHMWW}.

\begin{theorem}\label{Theorem(2,4)}
The following conditions are equivalent for any vector of nonnegative real numbers $(p_{12},p_{13},p_{14},p_{23},p_{24},p_{34})$.
\begin{enumerate}[(1)]\itemsep 5pt
\item There is a convex body $A\subseteq \mathbb{R}^4$ that satisfies
\[
p_{ij}=\big(\text{the area of the projection $\pi_{ij}(A)$}\big) \ \text{for all $i< j$.}
\]
\item There is a Euclidean triangle with side lengths $\sqrt{p_{12}p_{34}}, \sqrt{p_{13}p_{24}}, \sqrt{p_{14}p_{23}}$.
\end{enumerate}
\end{theorem}

In other words, $(p_{12},p_{13},p_{14},p_{23},p_{24},p_{34})$ is realizable as the vector of projection areas of a convex body in $\mathbb{R}^4$ if and only if it satisfies the triangle inequalities
\begin{align*}
\sqrt{p_{12}p_{34}} &\le  \sqrt{p_{13}p_{24}} + \sqrt{p_{14}p_{23}} \ \ \text{and}\\
\sqrt{p_{13}p_{24}} &\le  \sqrt{p_{12}p_{34}}  +  \sqrt{p_{14}p_{23}}  \ \ \text{and}\\
\sqrt{p_{14}p_{23}} &\le  \sqrt{p_{12}p_{34}}  + \sqrt{p_{13}p_{24}}. 
\end{align*}
A triangle is  said to be \emph{nondegenerate} when all the triangle inequalities are strict.

\begin{theorem}\label{Theorem(2,4)'}
The following conditions are equivalent for any vector of nonnegative real numbers $(p_{12},p_{13},p_{14},p_{23},p_{24},p_{34})$.
\begin{enumerate}[(1)]\itemsep 5pt
\item There is a smooth convex body $A\subseteq \mathbb{R}^4$ that satisfies
\[
p_{ij}=\big(\text{the area of the projection $\pi_{ij}(A)$}\big) \ \text{for all $i< j$.}
\]
\item There is a nondegenerate triangle with side lengths $\sqrt{p_{12}p_{34}}, \sqrt{p_{13}p_{24}}, \sqrt{p_{14}p_{23}}$.
\end{enumerate}
\end{theorem}

The \emph{realization space} of $(p_{12},p_{13},p_{14},p_{23},p_{24},p_{34})$ is the set of all convex bodies in $\mathbb{R}^4$ with the given projection areas. Theorem~\ref{Theorem(2,4)} characterizes the vectors with nonempty realization spaces.
One can show that, in a precise sense, the realization space of a given vector depends only on the associated triangle with side lengths $\sqrt{p_{12}p_{34}}, \sqrt{p_{13}p_{24}}, \sqrt{p_{14}p_{23}}$. 

\begin{example}
The vector $(\pi,\pi,\pi,\pi,\pi,\pi)$ is realizable as the vector of projection areas of a convex body in $\mathbb{R}^4$. For example, one may take the unit ball, the hypercube with side lengths $\sqrt{\pi}$, or more generally any convex body preserved by the $S_4$-symmetry of $\mathbb{R}^4$, scaled appropriately.
\end{example}

\begin{example}
According to Theorem~\ref{Theorem(2,4)}, the vector 
$(2,1,1,1,1,2)$
 is realizable as the vector of projection areas of a convex body $A \subseteq\mathbb{R}^4$ because there is a triangle with side lengths $2,1,1$. By Theorem~\ref{Theorem(2,4)'},  such a convex body  cannot be smooth. As a realization, one may take 
\[
A=\sqrt{2}\ \big(\text{the convex hull of $\mathbf{e}_1,\mathbf{e}_2,\mathbf{e}_3,\mathbf{e}_4,\mathbf{e}_1+\mathbf{e}_2,\mathbf{e}_3+\mathbf{e}_4$ in $\mathbb{R}^4$}\big).
\]
The projections $\pi_{12}(A)$ and $\pi_{34}(A)$ are squares with side lengths $\sqrt{2}$, and the remaining projections of $A$  
are triangles with side lengths $2,\sqrt{2},\sqrt{2}$.
\end{example}

\begin{example}
According to Theorem~\ref{Theorem(2,4)}, the vector $(3,2,1,1,2,3)$ is realizable as the vector of projection areas of a convex body $A \subseteq \mathbb{R}^4$ because there is a triangle with side lengths $3,2,1$. For example, one may take
the $4 \times 16$ matrices
\[
{\setlength{\arraycolsep}{3pt}
\begin{array}{r@{\;}c@{\;}l} 
L \coloneqq  \left[ \begin{array}{*{20}r} \phantom{-}0 & -1 & 0 & -1 & -1 & 0 & 0 & \phantom{-}0 & \phantom{-}1 & \phantom{-}1 & 1 & 0 & 0 & -1 & 1 & 0 \\ 0 & -1 & 0 & 1 & 1 & 0 & 0 & 0 & 1 & 1 & -1 & 0 & 0 & -1 & -1 & 0 \\ 1 & 0 & -1 & 0 & 0 & -1 & -1 & 1 & 0 & 0 & 0 & 1 & -1 & 0 & 0 & 1 \\ 1 & 0 & -1 & 0 & 0 & 1 & 1 & 1 & 0 & 0 & 0 & -1 & -1 & 0 & 0 & -1 \end{array} \right], \\[3em] 
M \coloneqq \left[ \begin{array}{*{20}r} 0 & \phantom{-} 0 & \phantom{-}0 & 0 & 1 & 0 & -1 & -1 & -1 & 0 & 0 & 1 & 1 & 1 & -1 & 0 \\ -2 & 0 & 0 & -2 & -1 & 0 & -1 & -1 & -1 & -2 & 0 & -1 & -1 & -1 & -1 & -2 \\ 0 & 0 & 0 & 0 & -1 & 0 & 1 & -1 & -1 & 0 & 0 & -1 & 1 & 1 & 1 & 0 \\ -2 & 0 & 0 & 0 & -1 & -2 & -1 & -1 & -1 & -2 & -2 & -1 & -1 & -1 & -1 & 0 \end{array} \right], \end{array}}
\]
and let $A$ be the convex hull in $\mathbb{R}^4$ of the $16$ columns of $L+\frac{1}{\sqrt{2}} M$.
With patience, one can check that $A$ has the  projections with given areas. For example, 
\[
\pi_{24}(A)=\big(\text{the convex hull of $\mathbf{e}_1,-\mathbf{e}_1,\mathbf{e}_3,-\mathbf{e}_3$ in $\mathbb{R}^2$}\big).
\]
Can you think of a simpler convex body in $\mathbb{R}^4$ that has the same six projection areas? See \cite[Section 4]{HHMWW} for a discussion of this particular case.
\end{example}

The condition characterizing the realizability of $(p_{12},p_{13},p_{14},p_{23},p_{24},p_{34})$ in Theorem~\ref{Theorem(2,4)}   is precisely the validity of the  \emph{Pl\"ucker relation} 
\[
p_{12}p_{34}-p_{13}p_{24}+p_{14}p_{23}=0 \ \ \text{for the Grassmannian $\textrm{Gr}(2,4)$},
\]
interpreted over the  triangular hyperfield $\mathbb{T}_2$ on the set of nonnegative real numbers 
\cite{BHKL-triangular}.
As noted in  \cite[Proposition 3.1]{HHMWW}, the condition is also equivalent to the statement that the nonnegative symmetric matrix
\[  
\begin{pmatrix}
0&p_{12}&p_{13}&p_{14}\\
p_{12}&0&p_{23}&p_{24}\\
p_{13}&p_{23}&0&p_{34}\\
p_{14}&p_{24}&p_{34}&0
\end{pmatrix} 
\  \text{has at most one positive eigenvalue.}
\]
Theorem~\ref{Theorem(2,4)} is  related to the fact that the set of quadratic volume polynomials $\mathbb{V}^2_n(\mathbb{R},k)$ is equal to the set of Lorentzian polynomials (Section~\ref{SectionProjective}).

\subsection{}

One quickly encounters interesting questions when trying to formulate similar realization problems in higher dimensions.
The asymmetry in the following pair of conjectures points to the distinction between \emph{volume polynomials} and \emph{covolume polynomials}  (Section~\ref{SectionVolumeCovolume}).

Let $\pi_{ij}$ be the coordinate projection of $\mathbb{R}^5$ onto the coordinate subspace orthogonal to the standard basis vectors $\mathbf{e}_i$ and $\mathbf{e}_j$.

\begin{conjecture}\label{Conjecture(2,5)}
The following conditions are equivalent for any vector of nonnegative real numbers $(p_{ij})_{1 \le i < j \le 5}$:\footnote{After the completion of this manuscript, Helen W. J. Zhang and Chengdong Zhao found an argument showing that there is no convex body \(A \subseteq \mathbb{R}^5\) satisfying
$p_{12}=p_{13}=p_{14}=p_{23}=p_{24}=p_{34}=1$ and 
$p_{15}=p_{25}=p_{35}=p_{45}=0$,
even though the corresponding matrix \(P(0)\) has eigenvalues \(3,0,-1,-1,-1\). On the other hand, Shouda Wang observed that there exist convex bodies \(A(t)\subseteq \mathbb{R}^5\) for \(t>0\) whose projection-volume matrices \(P(t)\) converge to the matrix \(P(0)\) as \(t\to 0\).
This leaves open the following variant of Conjecture~\ref{Conjecture(2,5)}: A $5\times 5$ nonnegative symmetric matrix with zero diagonal entries has at most one positive eigenvalue if and only if it is a limit of matrices of the form
$
(\operatorname{vol}\,\pi_{ij}(A))$, where \(A\) is a convex body in \(\mathbb{R}^5\).
}
\begin{enumerate}[(1)]\itemsep 5pt
\item There is a convex body $A\subseteq \mathbb{R}^5$ that satisfies
\[
p_{ij}=\big(\text{the volume of the projection $\pi_{ij}(A)$}\big) \ \text{for all $i< j$.}
\]
\item The nonnegative symmetric matrix
\[  
P \coloneqq \begin{pmatrix}
0&p_{12}&p_{13}&p_{14} & p_{15}\\
p_{12}&0&p_{23}&p_{24} & p_{25}\\
p_{13}&p_{23}&0&p_{34} & p_{35}\\
p_{14}&p_{24}&p_{34}&0 & p_{45}\\
p_{15}& p_{25}& p_{35}& p_{45}& 0
\end{pmatrix}
\ \text{has at most one positive eigenvalue.}
\]
\end{enumerate}
\end{conjecture}

Let $\varphi_{ij}$ be the coordinate projection of $\mathbb{R}^5$ onto the coordinate subspace spanned by the standard basis vectors $\mathbf{e}_i$ and $\mathbf{e}_j$. 

\begin{conjecture}\label{Conjecture(3,5)}
The following conditions are equivalent for any vector of nonnegative real numbers $(q_{ij})_{1 \le i < j \le 5}$:
\begin{enumerate}[(1)]\itemsep 5pt
\item There is a convex body $A\subseteq \mathbb{R}^5$ that satisfies
\[
q_{ij}=\big(\text{the area of the projection $\varphi_{ij}(A)$}\big) \ \text{for all $i< j$.}
\]
\item Every $4 \times 4$ principal submatrix of the nonnegative symmetric matrix
\[  
Q \coloneqq
\begin{pmatrix}
0&q_{12}&q_{13}&q_{14} & q_{15}\\
q_{12}&0&q_{23}&q_{24} & q_{25}\\
q_{13}&q_{23}&0&q_{34} & q_{35}\\
q_{14}&q_{24}&q_{34}&0 & q_{45}\\
q_{15}& q_{25}& q_{35}& q_{45}& 0
\end{pmatrix}
\ \text{has at most one positive eigenvalue.}
\]
\end{enumerate}
\end{conjecture}

In both cases, the forward implication follows from the Alexandrov--Fenchel inequality on mixed volumes. 
Similar conjectures can be made more generally for projections $\pi_{ij}:\mathbb{R}^d \to \mathbb{R}^{d-2}$
and $\varphi_{ij}:\mathbb{R}^d \to \mathbb{R}^2$.

\begin{example}\label{ExampleVolumeCovolume}
To compare Conjectures~\ref{Conjecture(2,5)} and ~\ref{Conjecture(3,5)}, we consider the case of $(4,1,1,1,1,1,1,1,1,1)$. Since  the corresponding symmetric matrix has eigenvalues
\[
 3+\sqrt{7}, \ 3-\sqrt{7}, \ -1, \ -1, \ -4,
 \]
there is no convex body $A$ in $\mathbb{R}^5$ such that the projection $\pi_{12}(A)$ has volume $4$ while all the other projections $\pi_{ij}(A)$ have volume $1$.
On the other hand, every $4 \times 4$ principal submatrix of the same matrix has at most one positive eigenvalue, suggesting that there is a convex body $A$ in $\mathbb{R}^5$ such that the projection $\varphi_{12}(A)$ has area $4$ while all the other projections $\varphi_{ij}(A)$ have area $1$. Indeed, one can take
\[
A=\big(\text{the convex hull of $2\mathbf{e}_1,2\mathbf{e}_2,2\mathbf{e}_1+2\mathbf{e}_2,\mathbf{e}_3,\mathbf{e}_4,\mathbf{e}_5,\mathbf{e}_3+\mathbf{e}_4,\mathbf{e}_3+\mathbf{e}_5,\mathbf{e}_4+\mathbf{e}_5$}\big),
\]
which is consistent with Conjectures~\ref{Conjecture(2,5)} and ~\ref{Conjecture(3,5)}.
\end{example}

\subsection{}

We may pose analogous realization problems in the setting of projective geometry.
Fix an algebraically closed field $k$, and consider the projective line $\mathbb{P}^1$ over $k$.
If $S$ is an irreducible surface in $(\mathbb{P}^1)^4$, we can uniquely express its homology class as a nonnegative integral linear combination
\[
[S]=p_{12}[\PP^1 \times \PP^1 \times \PP^0 \times \PP^0]
+\cdots +p_{34}[\PP^0 \times \PP^0 \times \PP^1 \times \PP^1] \in \textrm{CH}((\mathbb{P}^1)^4).
\] 
Which vectors of nonnegative integers $(p_{12},p_{13},p_{14},p_{23},p_{24},p_{34})$ can arise in this way?
We call such homology classes \emph{realizable}.\footnote{Such questions are algebraic analogues of the \emph{Steenrod problem} in topology \cite[Problem 25]{Eilenberg}, which asks whether every homology class in any simplicial complex $X$ is the image of the fundamental class of a closed oriented manifold by a map into the simplicial complex.} 
While the answer to this question is not known, the following partial result was given in \cite[Theorem 1.6]{HHMWW}.

\begin{theorem}\label{Theorem(2,4)''}
The following conditions are equivalent for any vector of nonnegative rational numbers $(p_{12},p_{13},p_{14},p_{23},p_{24},p_{34})$.
\begin{enumerate}[(1)]\itemsep 5pt
\item There is an irreducible surface $S \subseteq (\mathbb{P}^1)^4$ and  a nonnegative $\lambda \in \mathbb{Q}$ such that
\[
\lambda [S]=p_{12}[\PP^1 \times \PP^1 \times \PP^0 \times \PP^0]
+\cdots +p_{34}[\PP^0 \times \PP^0 \times \PP^1 \times \PP^1].
\] 
\item There is a Euclidean triangle with side lengths $\sqrt{p_{12}p_{34}}, \sqrt{p_{13}p_{24}}, \sqrt{p_{14}p_{23}}$.
\end{enumerate}
\end{theorem}

The algebraic realization problem here has additional obstructions not present in the convex realization problem in Section~\ref{SubsectionConvex}.
For example, the homology class corresponding to $(1,1,1,1,1,3)$ is not realizable by an irreducible surface in $(\mathbb{P}^1)^4$, although there is a Euclidean triangle with side lengths $\sqrt{3},1,1$. 
To see this, note that the hypothetical surface $S$ should satisfy
\[
[\pi_{3}(S)]=[\pi_{4}(S)]=[\mathbb{P}^1 \times \mathbb{P}^1 \times \mathbb{P}^0]+[\mathbb{P}^1 \times \mathbb{P}^0 \times \mathbb{P}^1]+[\mathbb{P}^0 \times \mathbb{P}^1 \times \mathbb{P}^1],
\]
where $\pi_i$ is the projection $(\mathbb{P}^1)^4 \to (\mathbb{P}^1)^3$ that forgets the $i$-th coordinate of $(\mathbb{P}^1)^4$. 
Thus, the defining equations of the hypersurfaces  $\pi_{3}^{-1} \pi_{3}(S)$ and $\pi_{4}^{-1} \pi_{4}(S)$ in an affine chart of $(\mathbb{P}^1)^4$
are of the form 
\begin{align*}
*1 + *x_1+  *x_2+ *x_4+ *x_1x_2+ *x_1x_4+ *x_2x_4+*x_1x_2x_4 &=0,\\
   *1 + *x_1+  *x_2+ *x_3+ *x_1x_2+ *x_1x_3+ *x_2x_3+*x_1x_2x_3 &=0, 
\end{align*}
where the $*$'s are placeholders for the coefficients in $k$. 
For generic values of $x_3$ and $x_4$, the displayed system of equations reduces to
\begin{align*}
*1 + *x_1+  *x_2+*x_1x_2&=0,\\
   *1 + *x_1+  *x_2+*x_1x_2 &=0, 
\end{align*}
which has at most 2 solutions. This contradicts the assumption that  $p_{34}=3$.
 The proof of Theorem~\ref{Theorem(2,4)''} in \cite{HHMWW} shows that, in fact, the homology class corresponding to $(2,2,2,2,2,6)$ is realizable by an irreducible surface in $(\mathbb{P}^1)^4$.

Modulo this subtlety involving integral coefficients, one can formulate statements in algebraic geometry parallel to Conjectures \ref{Conjecture(2,5)} and \ref{Conjecture(3,5)} in convex geometry, the first of which is a special case of  \cite[Theorem 1.8]{HHMWW}.

\begin{theorem}
The following conditions are equivalent for any vector of nonnegative integers $(p_{ij})_{1 \le i < j \le 5}$:
\begin{enumerate}[(1)]\itemsep 5pt
\item There is an irreducible surface  $S \subseteq (\mathbb{P}^1)^5$ and  a nonnegative $\lambda \in \mathbb{Q}$ such that
\[
\lambda [S]=p_{12}[\PP^1 \times \PP^1 \times \PP^0 \times \PP^0 \times \PP^0]
+\cdots +p_{45}[\PP^0 \times \PP^0\times \PP^0 \times \PP^1 \times \PP^1].
\] 
\item The nonnegative symmetric matrix
\[  
\begin{pmatrix}
0&p_{12}&p_{13}&p_{14} & p_{15}\\
p_{12}&0&p_{23}&p_{24} & p_{25}\\
p_{13}&p_{23}&0&p_{34} & p_{35}\\
p_{14}&p_{24}&p_{34}&0 & p_{45}\\
p_{15}& p_{25}& p_{35}& p_{45}& 0
\end{pmatrix}
\ \text{has at most one positive eigenvalue.}
\]
\end{enumerate}
\end{theorem}

The following conjecture suggests the possibility that, in Conjecture~\ref{Conjecture(3,5)}, the convex body $A$ can be chosen to be a rational convex polytope whenever all $q_{ij}$ are rational.

\begin{conjecture}\label{ConjectureCovolume}
The following conditions are equivalent for any vector of nonnegative integers $(q_{ij})_{1 \le i < j \le 5}$:
\begin{enumerate}[(1)]\itemsep 5pt
\item There is an irreducible threefold  $S \subseteq (\mathbb{P}^1)^5$ and  a nonnegative $\lambda \in \mathbb{Q}$ such that
\[
\lambda [S]=q_{12}[\PP^0 \times \PP^0 \times  \PP^1 \times \PP^1 \times \PP^1]
+\cdots +q_{45}[\PP^1 \times \PP^1\times \PP^1 \times \PP^0 \times \PP^0].
\] 
\item Every $4 \times 4$ principal submatrix of the nonnegative symmetric matrix
\[  
\begin{pmatrix}
0&q_{12}&q_{13}&q_{14} & q_{15}\\
q_{12}&0&q_{23}&q_{24} & q_{25}\\
q_{13}&q_{23}&0&q_{34} & q_{35}\\
q_{14}&q_{24}&q_{34}&0 & q_{45}\\
q_{15}& q_{25}& q_{35}& q_{45}& 0
\end{pmatrix}
\ \text{has at most one positive eigenvalue.}
\]
\end{enumerate}
\end{conjecture}

In both cases, the Hodge index theorem for projective surfaces can be used to show the forward implication.

\section{Volume polynomials in convex geometry}\label{SectionConvex}

\subsection{}

In  \cite{Minkowski}, Minkowski made a foundational observation that has since become a cornerstone of convex geometry: 
\begin{quote}
\emph{The volume of the Minkowski sum of convex bodies varies polynomially under scaling.}
\end{quote}
More precisely,  for positive integers $n$ and $d$, and any collection of $n$ convex bodies $C=(C_1,\ldots,C_n)$ in the $d$-dimensional Euclidean space $\RR^d$, the function
\[
f_{C}: \RR^n_{\ge 0} \longrightarrow \RR_{\ge 0}, \qquad (x_1,\ldots,x_n) \longmapsto \frac{1}{d!} \vol(x_1C_1+\dots+x_nC_n)
\]
is a degree $d$ homogeneous polynomial in $x=(x_1,\ldots,x_n)$.
This polynomial, called the \emph{volume polynomial} of $C$, is then used to define the \emph{mixed volume} of convex bodies in $C$ as its normalized coefficients
\[
\textrm{MV}(C_{i_1},\ldots,C_{i_d}) \coloneqq  \frac{\partial}{\partial x_{i_1}}\cdots \frac{\partial}{\partial x_{i_d}} f_{C}(x_1,\ldots,x_n). 
\]
The mixed volume is symmetric in its arguments, and it is multilinear with respect to the Minkowski sum and nonnegative scaling: For any convex bodies $B_1$ and $B_2$ in $\mathbb{R}^d$ and nonnegative real numbers $\lambda_1$ and  $\lambda_2$, we have
\[
\textrm{MV}(\lambda_1 B_1+\lambda_2 B_2,C_{2},\ldots,C_{d})=\lambda_1 \textrm{MV}(B_1,C_{2},\ldots,C_{d})+\lambda_2 \textrm{MV}(B_2,C_{2},\ldots,C_{d}).
\]
When all the arguments coincide, the mixed volume reduces to the usual volume: For any convex body $A$ in $\mathbb{R}^d$, we have
\[
\textrm{MV}(\underbrace{A,\dots,A}_{d}) = \vol(A).
\]
More generally, 
if $I_j$ is the unit interval joining the origin and  $\mathbf{e}_j$ in $\mathbb{R}^d$, we have
\[
\binom{d}{k} \, \textrm{MV}(I_1,\ldots, I_k,\underbrace{A,\dots,A}_{d-k}) = \frac{1}{k!} \vol( \pi_{1 \cdots k} A),
\]
where $\pi_{1 \cdots k}$ is the projection onto the coordinate subspace orthogonal to $\mathbf{e}_1,\ldots,\mathbf{e}_k$. 
For a comprehensive introduction to mixed volumes, see \cite[Chapter 5]{Schneider}.

Mixed volumes of convex bodies satisfy a rich collection of fundamental inequalities, the most basic being \emph{nonnegativity}:
\[
0 \le  \textrm{MV}(C_1,\ldots,C_d) \ \text{for any convex bodies $C_1,\ldots,C_d$ in $\mathbb{R}^d$.}
\]
More generally, mixed volumes are \emph{monotone} in each argument: If $B_i \subseteq C_i$ are convex bodies in $\mathbb{R}^d$, we have
\[
 \textrm{MV}(B_1,\ldots,B_d)  \le  \textrm{MV}(C_1,\ldots,C_d). 
\]
Apart from the nonnegativity, 
the most important inequality involving mixed volumes is the \emph{Alexandrov--Fenchel inequality}, which generalizes classical inequalities such as the isoperimetric and Brunn--Minkowski inequalities. It states that, for any convex bodies $C_1,\ldots,C_d$ in $\mathbb{R}^d$, we have
\[
 \textrm{MV}(C_1,C_1,C_3,\ldots,C_d) \, \textrm{MV}(C_2,C_2,C_3,\ldots,C_d) \le  \textrm{MV}(C_1,C_2,C_3,\ldots,C_d)^2.
\]
This inequality is a cornerstone of modern convex geometry and underlies many structural and analytic results.

\subsection{}

The \emph{realization problem} for volume polynomials of convex bodies is to determine which homogeneous polynomials of degree $d$ in $n$ variables with nonnegative coefficients can be realized as the volume polynomials of $n$ convex bodies in $\mathbb{R}^d$. The problem of finding the full set of inequalities for mixed volumes is sometimes referred to as \emph{Alexandrov's problem}.

For a degree $d$ homogeneous polynomial $f$ in $n$ variables $x=(x_1,\ldots,x_n)$, we write
\[
f(x)=\sum_{\alpha \in \Delta^d_n} p_\alpha x^{[\alpha]}, \quad  x^{[\alpha]} \coloneqq \frac{x^\alpha}{\alpha!}= \frac{x_1^{\alpha_1}}{\alpha_1!} \cdots  \frac{x_n^{\alpha_n}}{\alpha_n!},
\]
where $p_\alpha$ are the normalized coefficients of $f$ and $\Delta^d_n$ is the discrete simplex consisting of all the nonnegative vectors in $\mathbb{Z}^n$ whose coordinates sum to $d$.
Here are the first two necessary conditions for $f$ to be a volume polynomial.
\begin{enumerate}[(1)]\itemsep 5pt
\item (Nonnegative change of coordinates) If $f(x)$ is a volume polynomial of $n$ convex bodies, then $f(Ay)$ is a volume polynomial of $m$ convex bodies, for any $n \times m$ nonnegative matrix $A$ and variables $y=(y_1,\ldots,y_m)$.
\item (Alexandrov--Fenchel inequality) If $f(x)$ is a volume polynomial of convex bodies, then its normalized coefficients satisfy 
\[
p_{\alpha+\mathrm{e}_i-\mathrm{e}_j}p_{\alpha-\mathrm{e}_i+\mathrm{e}_j } \le p_\alpha^2 \ \ \text{for any $\alpha \in \Delta^d_n$ and any $1 \le i<j \le n$.}
\]
\end{enumerate}
The first condition is a formal consequence of the observation that the Minkowski sum of convex bodies is a convex body. The combination of the two properties leads to the conclusion that any volume polynomial must be a \emph{Lorentzian polynomial} in the sense of \cite{BrandenHuh}. 
Here we give an equivalent definition, following \cite[Section 2]{BrandenLeake}.
We set
\begin{multline*}
\big\{\text{Lorentzian polynomials of degree $\le 1$}\big\}=\\
\big\{\text{homogeneous polynomials of degree $\le 1$ with nonnegative coefficients}\big\}.
\end{multline*}
For a nonnegative vector $u=(u_1,\ldots,u_n)$, we write $\partial_u$ for the corresponding directional derivative $\sum_{i=1}^n u_i  \partial_i$.

\begin{definition}\label{DefinitionLorentzian}
A homogeneous polynomial $f$ of degree $d \ge 2$ in $n$ variables with nonnegative coefficients is \emph{Lorentzian} if, for all  $v_1,\ldots,v_d \in \mathbb{R}^n_{\ge 0}$, we have
\[
 \big( \partial_{v_1}  \partial_{v_1} \partial_{v_3}  \cdots \partial_{v_d} f \big)  \big( \partial_{v_2}  \partial_{v_2} \partial_{v_3}  \cdots \partial_{v_d} f\big)   \le \big( \partial_{v_1} \partial_{v_2} \partial_{v_3}\cdots \partial_{v_d} f \big)^2.
\]
\end{definition}

Applying the Alexandrov--Fenchel inequality after a nonnegative linear change of coordinates, we see that
\[
\big\{\text{volume polynomials of convex bodies}\big\} \subseteq \big\{\text{Lorentzian polynomials}\big\}.
\]
Thus, Alexandrov's problem is to find inequalities between mixed volumes that identify the volume polynomials of convex bodies among Lorentzian polynomials.

\begin{example}[$n=2$]\label{ExampleBivariate}
According to  \cite[Example 2.26]{BrandenHuh}, a bivariate polynomial with nonnegative coefficients 
\[
f= \sum_{a=0}^d p_a\, \frac{x_1^a}{a!} \frac{x_2^{d-a}}{(d-a)!}
\]
is Lorentzian if and only if the sequence $p_0,\ldots,p_d$ has no internal zeros and 
\[
p_{a-1}p_{a+1} \le p_a^2 \ \ \text{for all positive integers $a<d$}.
\]
In  \cite{Shephard}, Shephard showed that any such polynomial is the volume polynomial of two convex bodies in $\mathbb{R}^d$.
This characterizes volume polynomials of convex bodies in two variables:
\begin{quote}
\emph{A homogeneous polynomial in two variables is the volume polynomial of two convex bodies if and only if it is Lorentzian.}
\end{quote} 
When every  $p_i$ is rational, Shephard's construction gives two rational convex polytopes. 
This is used in \cite[Theorem 21]{HuhChromatic} to characterize realizable homology classes in $\mathbb{P}^d \times \mathbb{P}^d$ up to a multiple:
Some nonnegative rational multiple of the class
\[
\sum_{a=0}^d p_a \, [\mathbb{P}^a \times \mathbb{P}^{d-a}]   \in \textrm{CH}(\mathbb{P}^d \times \mathbb{P}^d)
\]
is the class of an irreducible subvariety if and only if  $p_0,\ldots,p_d$ is a log-concave sequence of nonnegative rational numbers with no internal zeros.\footnote{As observed in \cite[Section 5]{HuhCorrespondence}, there is no irreducible subvariety of $\mathbb{P}^5 \times \mathbb{P}^5$ whose homology class corresponds to the log-concave sequence $(1,2,3,4,2,1)$.}
\end{example}

\begin{example}[$d=1$]
By definition, a linear form is Lorentzian if and only if all its coefficients are nonnegative. 
Any such linear form is a volume polynomial of convex bodies in $\mathbb{R}^1$:
\[
\vol(x_1 C_1+\cdots+x_nC_n)=x_1 \vol(C_1)+\cdots+x_n\vol(C_n), \ \ C_1,\ldots,C_n \subseteq \mathbb{R}^1.
\]
\end{example}

\begin{example}[$d=2$]\label{ExampleQuadratic}
A quadratic form 
is Lorentzian if and only if all its coefficients are nonnegative and its Hessian has at most one positive eigenvalue \cite[Section 2]{BrandenHuh}.
In \cite{Heine}, Heine showed that, when there are at most three variables, any such quadratic form is the volume polynomial of three convex bodies in $\mathbb{R}^2$.
This characterizes quadratic volume polynomials of convex bodies in three variables:
\begin{quote}
\emph{A ternary quadratic form is the volume polynomial of three convex bodies if and only if it is Lorentzian.}
\end{quote}
The analogous statement fails when  $n=4$. For example, as observed in \cite[Theorem 5]{Shephard}, there are no convex bodies $C_1,C_2,C_3,C_4$ in $\mathbb{R}^2$ satisfying  
\[
\vol(x_1C_1+x_2C_2+x_3C_3+x_4C_4)=x_1x_2+x_1x_3+x_1x_4+x_2x_3+x_2x_4+x_3x_4,
\]
even though the right-hand side is a Lorentzian polynomial. 
 In fact,  using the compactness theorem of Shephard for the affine equivalence classes of convex bodies \cite[Theorem 1]{Shephard}, one can show that the displayed elementary symmetric polynomial is not even the limit of volume polynomials of convex bodies in the plane.
This contrasts with the fact that there is an irreducible surface in $(\mathbb{P}^1)^4$ with class $(1,1,1,1,1,1)$ in the Chow group. For example, one may take the closure of a general two-dimensional linear subspace of an affine chart of $(\mathbb{P}^1)^4$.
\end{example}

\subsection{}\label{SectionLorentzian}

The main result of \cite{BrandenHuh} provides a finite description of the set of Lorentzian polynomials that generalizes Example~\ref{ExampleBivariate}.
The central notion is that of a generalized permutohedron.
Let $E$ be a finite set with $n$ elements, and let $\{e_i\}_{i \in E}$ be the standard basis of $\mathbb{R}^E$.

\begin{definition}
A \emph{generalized permutohedron} is a polytope in $\mathbb{R}^E$ all of whose edges are in the direction $e_i-e_j$ for some $i$ and $j$ in $E$.
\end{definition}

A generalized permutohedron is \emph{integral} if all its vertices belong to $\mathbb{Z}^E \subseteq \mathbb{R}^E$. 
For example, the \emph{standard permutohedron} in $\mathbb{R}^n$, which is the convex hull of all permutations of $(1,2,\ldots,n)$,
and the \emph{$k$-th hypersimplex} in $\mathbb{R}^n$, which is the convex hull of all permutations of $(\underbrace{1,\ldots,1}_{k},\underbrace{0,\ldots,0}_{n-k})$, are integral generalized permutohedra.
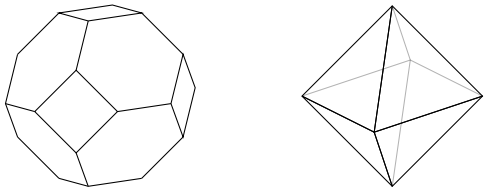
\[
\begin{tikzpicture}[z={(-.1,0,.5)},scale=0.55]
\foreach [var=\x, var=\y, var=\z, count=\n] in {
2/1/0,2/0/1,2/-1/0,2/0/-1,
1/0/2,0/1/2,-1/0/2,0/-1/2,
1/2/0,0/2/1,-1/2/0,0/2/-1,
-2/1/0,-2/0/1,-2/-1/0,-2/0/-1,
1/0/-2,0/1/-2,-1/0/-2,0/-1/-2,
1/-2/0,0/-2/1,-1/-2/0,0/-2/-1
}{\coordinate (n\n) at (\x,\y,\z);}
\draw (n1)--(n2)--(n3)--(n4)--cycle;
\draw (n5)--(n6)--(n7)--(n8)--cycle;
\draw (n9)--(n10)--(n11)--(n12)--cycle;
\draw (n21)--(n22)--(n23);
\draw (n13)--(n14)--(n15);
\draw (n6)--(n10);
\draw (n2)--(n5);
\draw (n8)--(n22);
\draw (n15)--(n23);
\draw (n7)--(n14);
\draw (n11)--(n13);
\draw (n1)--(n9);
\draw (n3)--(n21);
\end{tikzpicture}
\qquad \qquad
\begin{tikzpicture}[scale=2.4]
    \coordinate (A1) at (0,0);
    \coordinate (A2) at (0.6,0.2);
    \coordinate (A3) at (1,0);
    \coordinate (A4) at (0.4,-0.2);
    \coordinate (B1) at (0.5,0.5);
    \coordinate (B2) at (0.5,-0.5);

   \begin{scope}[opacity=0.3]
        \draw (A1) -- (A2) -- (A3);
        \draw (B1) -- (A2) -- (B2);
    \end{scope}

    \draw (A1) -- (A4) -- (B1);
    \draw (A1) -- (A4) -- (B2);
    \draw (A3) -- (A4) -- (B1);
    \draw (A3) -- (A4) -- (B2);
    \draw (B1) -- (A1) -- (B2) -- (A3) --cycle;
\end{tikzpicture}
\]
The above pictures show the standard permutohedron and the second hypersimplex in $\mathbb{R}^4$.
Generalized permutohedra are precisely the 
polytopes obtained from the standard permutohedron by moving the vertices so that all the edge directions are preserved \cite{Postnikov}.

\begin{definition}\label{DefinitionM}
A subset $J \subseteq \mathbb{Z}^E_{\ge 0}$ is \emph{$\mathrm{M}$-convex} if it is the set of all lattice points of an integral generalized permutohedron.
A \emph{matroid} on $E$ is an $\mathrm{M}$-convex subset of $\mathbb{Z}^E_{\ge 0}$ consisting of zero-one vectors.
The vectors in a matroid $J$ are called \emph{bases} of $J$.
\end{definition}

The notion of $\mathrm{M}$-convex sets originates in discrete convex analysis \cite{Murota}. 
In \cite[Chapter 4]{Murota}, one can find several other equivalent characterizations of $\mathrm{M}$-convex sets.
For example, a subset $J \subseteq \mathbb{Z}^E_{\ge 0}$ is $\mathrm{M}$-convex exactly when it satisfies the \emph{symmetric basis exchange property}: 
\begin{quote}
\emph{For any $\alpha, \beta \in J$ and $i \in E$ with $\alpha_i>\beta_i$, there is  $j \in E$ with
\[
\alpha_j<\beta_j \ \ \text{and} \ \ \alpha-e_i+e_j \in J \ \ \text{and} \ \ \beta-e_j+e_i \in J.
\]}
\end{quote}
 For background specific to  matroids, see \cite{Oxley}. 
 

\begin{definition}\label{def:polymatroid}
A function $h:2^E \to \ZZ_{\ge 0}$ is a \emph{polymatroid rank function} if it satisfies the following properties:
\begin{enumerate}[(1)]\itemsep 5pt
  \item \emph{Normalization:} $h(\varnothing) = 0$.
  \item \emph{Monotonicity:} 
    $h(A) \le h(B)$ for all $A \subseteq B \subseteq E$.
  \item \emph{Submodularity:} 
    $h(A\cup B) + h(A\cap B) \le h(A) + h(B)$ for all $A,B \subseteq E$.
\end{enumerate}
A polymatroid rank function $h$ is a \emph{matroid rank function} if $h(A) \le |A|$ for all $A\subseteq E$.
\end{definition}

We recall the standard bijection between polymatroid rank functions on $E$ and nonempty $\mathrm{M}$-convex subsets of $\ZZ^E_{\ge 0}$ from \cite[Chapter 4]{Murota}. 
For  $A \subseteq E$ and $\alpha \in \ZZ^E_{\ge 0}$, we set $\alpha_A \coloneqq \sum_{i \in A} \alpha_i$.
\begin{enumerate}[(1)]\itemsep 5pt
\item A polymatroid rank function $h$ defines 
\[
J_{h}
\coloneqq
\left\{
\alpha \in \mathbb{Z}_{\ge0}^E
\;\middle|\;
\alpha_E= h(E) \ \ \text{and} \ \  
\alpha_A \le h(A)\ \text{for all $A\subseteq E$}
\right\},
\]
which is an $\mathrm{M}$-convex subset of $\ZZ^E_{\ge 0}$.
\item  An $\mathrm{M}$-convex subset $J$ of $\ZZ^E_{\ge 0}$ defines 
\[
h_J:2^E \longrightarrow \ZZ_{\ge 0}, \quad h_J(A)\coloneqq \max\big\{ \beta_A  \mid \text{$\beta\le \alpha$ for some $\alpha \in J$}\big\},
\]
which is a polymatroid rank function on $E$.
\end{enumerate}

The constructions $J_h$ and $h_J$ are mutually inverse, providing a polymatroid generalization of the classical cryptomorphism between the  matroid rank function axioms and the symmetric basis exchange property.
A \emph{polymatroid} $\mathscr{P}$ is a pair $(h=h_J,J=J_h)$, where
 $h$ is the \emph{rank function} of $\mathscr{P}$ and $J$ is the \emph{set of bases} of $\mathscr{P}$. 
A polymatroid $\mathscr{P}$ is a \emph{matroid} if $h$ is a matroid rank function, or equivalently if $J$ consists of zero-one vectors.
Throughout this text, we restrict attention to integral polymatroids and do not consider nonintegral ones. Accordingly, we use the terms \emph{polymatroid} and 
\emph{$\mathrm{M}$-convex set} interchangeably.



\begin{example}[Graphic matroids]\label{ExampleGraphicMatroid}
For any finite connected graph $G$ with the edge set $E$, consider the set of indicator vectors
\[
J(G) \coloneqq \{ e_B \ | \ \text{$B$ is a spanning tree of $G$}\} \subseteq \mathbb{Z}^{E}_{\ge 0}.
\]
The subset $J(G)$ is $\mathrm{M}$-convex for any such $G$.
Such matroids are said to be \emph{graphic}.
\end{example}

\begin{example}[Linear matroids]\label{ExampleRepresentableMatroid}
For any function $\varphi:E \to W$ from a finite set $E$ to a vector space $W$ over a field $k$, consider the set of indicator vectors
\[
J(\varphi)\coloneqq \{e_B \ | \ \text{$\varphi(B)$ is a basis of $W$}\} \subseteq \mathbb{Z}^{E}_{\ge 0}.
\]
The subset $J(\varphi)$ is $\mathrm{M}$-convex for any $\varphi:E \to W$.
Such matroids are said to be \emph{linear over $k$},
and the function $\varphi$ is called a \emph{linear realization over $k$}.
One typically requires without loss of generality that the image of $\varphi$ spans $W$.
A graphic matroid is linearly  realizable over every field \cite[Section 5.1]{Oxley}.
In general, a matroid may or may not have a linear realization over $k$:
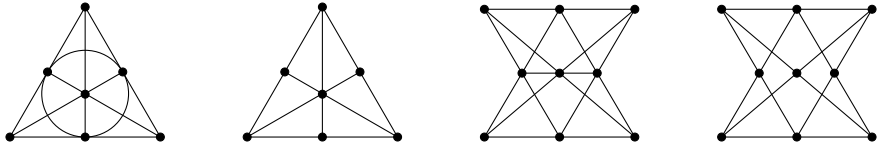
\begin{center}
\begin{tikzpicture}
\draw (0,0) -- (2,0);
\draw (2,0) -- (1,1.73);
\draw (1,1.73) -- (0,0);
\draw (0,0) -- (1.5,.866);
\draw (2,0) -- (.5,.866);
\draw (1,1.73) -- (1,0);
\draw (1,0.577) circle [radius=0.577];
\draw [fill=black] (0,0) circle (1.5pt); 
\draw [fill=black] (2,0) circle (1.5pt); 
\draw [fill=black] (1,1.73) circle (1.5pt); 
\draw [fill=black] (1.5,.866) circle (1.5pt); 
\draw [fill=black] (.5,.866) circle (1.5pt); 
\draw [fill=black] (1,0) circle (1.5pt); 
\draw [fill=black] (1,0.57) circle (1.5pt); 
\end{tikzpicture}
 \hspace{0.8cm}
\begin{tikzpicture}
\draw (0,0) -- (2,0);
\draw (2,0) -- (1,1.73);
\draw (1,1.73) -- (0,0);
\draw (0,0) -- (1.5,.866);
\draw (2,0) -- (.5,.866);
\draw (1,1.73) -- (1,0);
\draw [fill=black] (0,0) circle (1.5pt); 
\draw [fill=black] (2,0) circle (1.5pt); 
\draw [fill=black] (1,1.73) circle (1.5pt); 
\draw [fill=black] (1.5,.866) circle (1.5pt); 
\draw [fill=black] (.5,.866) circle (1.5pt); 
\draw [fill=black] (1,0) circle (1.5pt); 
\draw [fill=black] (1,0.57) circle (1.5pt); 
\end{tikzpicture}
\hspace{0.8cm}
\begin{tikzpicture}
\draw (0,0) -- (2,0);
\draw (0,1.7) -- (2,1.7);
\draw (0,0) -- (1,1.7);
\draw (0,0) -- (2,1.7);
\draw (1,0) -- (0,1.7);
\draw (1,0) -- (2,1.7);
\draw (2,0) -- (0,1.7);
\draw (2,0) -- (1,1.7);

\draw (0.5,.85) -- (1.5,.85);

\draw [fill=black] (0,0) circle (1.5pt); 
\draw [fill=black] (1,0) circle (1.5pt); 
\draw [fill=black] (2,0) circle (1.5pt); 
\draw [fill=black] (0.5,.85) circle (1.5pt); 
\draw [fill=black] (1,.85) circle (1.5pt); 
\draw [fill=black] (1.5,.85) circle (1.5pt); 
\draw [fill=black] (0,1.7) circle (1.5pt); 
\draw [fill=black] (1,1.7) circle (1.5pt); 
\draw [fill=black] (2,1.7) circle (1.5pt); 
\end{tikzpicture}
\hspace{0.8cm}
\begin{tikzpicture}
\draw (0,0) -- (2,0);
\draw (0,1.7) -- (2,1.7);
\draw (0,0) -- (1,1.7);
\draw (0,0) -- (2,1.7);
\draw (1,0) -- (0,1.7);
\draw (1,0) -- (2,1.7);
\draw (2,0) -- (0,1.7);
\draw (2,0) -- (1,1.7);
\draw [fill=black] (0,0) circle (1.5pt); 
\draw [fill=black] (1,0) circle (1.5pt); 
\draw [fill=black] (2,0) circle (1.5pt); 
\draw [fill=black] (0.5,.85) circle (1.5pt); 
\draw [fill=black] (1,.85) circle (1.5pt); 
\draw [fill=black] (1.5,.85) circle (1.5pt); 
\draw [fill=black] (0,1.7) circle (1.5pt); 
\draw [fill=black] (1,1.7) circle (1.5pt); 
\draw [fill=black] (2,1.7) circle (1.5pt); 
\end{tikzpicture}
\end{center}
Among the four matroids pictured above, where the bases are given by all triples of points not on a line, the first is linear over $k$ if and only if the characteristic of $k$ is $2$, the second is linear over $k$ if and only if the characteristic of $k$ is not $2$, the third is linear over $k$ if and only if the cardinality of $k$ is not $2$, $3$, or $5$, and the fourth is not linear over any field \cite[Appendix]{Oxley}.
\end{example}

\begin{example}[Algebraic matroids]\label{ExampleAlgebraicMatroids}
For any function $\varphi:E \to \ell$ from a finite set $E$ to a field extension $\ell$ of $k$, consider the set of indicator vectors
\[
J(\varphi)\coloneqq \{e_B \ | \ \text{$\varphi(B)$ is a transcendence basis of $\ell$ over $k$}\} \subseteq \mathbb{Z}^{E}_{\ge 0}.
\]
The subset $J(\varphi)$ is $\mathrm{M}$-convex for any $\varphi:E \to \ell$.
Such matroids are said to be \emph{algebraic over $k$},
and the function $\varphi$ is called an \emph{algebraic realization over $k$}.
One typically requires without loss of generality that the image of $\varphi$ contains a transcendence basis of $\ell$ over $k$.
A linear matroid over $k$ is algebraic over $k$ \cite[Section 6.7]{Oxley}.
In general, a matroid may or may not have an algebraic realization over $k$:
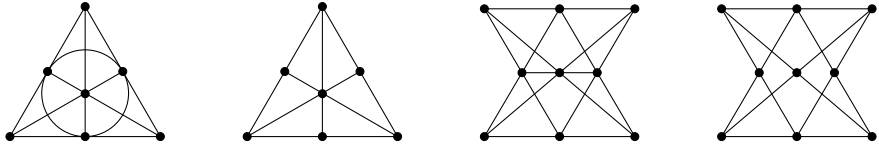
\begin{center}
\begin{tikzpicture}
\draw (0,0) -- (2,0);
\draw (2,0) -- (1,1.73);
\draw (1,1.73) -- (0,0);
\draw (0,0) -- (1.5,.866);
\draw (2,0) -- (.5,.866);
\draw (1,1.73) -- (1,0);
\draw (1,0.577) circle [radius=0.577];
\draw [fill=black] (0,0) circle (1.5pt); 
\draw [fill=black] (2,0) circle (1.5pt); 
\draw [fill=black] (1,1.73) circle (1.5pt); 
\draw [fill=black] (1.5,.866) circle (1.5pt); 
\draw [fill=black] (.5,.866) circle (1.5pt); 
\draw [fill=black] (1,0) circle (1.5pt); 
\draw [fill=black] (1,0.57) circle (1.5pt); 
\end{tikzpicture}
 \hspace{0.8cm}
\begin{tikzpicture}
\draw (0,0) -- (2,0);
\draw (2,0) -- (1,1.73);
\draw (1,1.73) -- (0,0);
\draw (0,0) -- (1.5,.866);
\draw (2,0) -- (.5,.866);
\draw (1,1.73) -- (1,0);
\draw [fill=black] (0,0) circle (1.5pt); 
\draw [fill=black] (2,0) circle (1.5pt); 
\draw [fill=black] (1,1.73) circle (1.5pt); 
\draw [fill=black] (1.5,.866) circle (1.5pt); 
\draw [fill=black] (.5,.866) circle (1.5pt); 
\draw [fill=black] (1,0) circle (1.5pt); 
\draw [fill=black] (1,0.57) circle (1.5pt); 
\end{tikzpicture}
\hspace{0.8cm}
\begin{tikzpicture}
\draw (0,0) -- (2,0);
\draw (0,1.7) -- (2,1.7);
\draw (0,0) -- (1,1.7);
\draw (0,0) -- (2,1.7);
\draw (1,0) -- (0,1.7);
\draw (1,0) -- (2,1.7);
\draw (2,0) -- (0,1.7);
\draw (2,0) -- (1,1.7);

\draw (0.5,.85) -- (1.5,.85);

\draw [fill=black] (0,0) circle (1.5pt); 
\draw [fill=black] (1,0) circle (1.5pt); 
\draw [fill=black] (2,0) circle (1.5pt); 
\draw [fill=black] (0.5,.85) circle (1.5pt); 
\draw [fill=black] (1,.85) circle (1.5pt); 
\draw [fill=black] (1.5,.85) circle (1.5pt); 
\draw [fill=black] (0,1.7) circle (1.5pt); 
\draw [fill=black] (1,1.7) circle (1.5pt); 
\draw [fill=black] (2,1.7) circle (1.5pt); 
\end{tikzpicture}
\hspace{0.8cm}
\begin{tikzpicture}
\draw (0,0) -- (2,0);
\draw (0,1.7) -- (2,1.7);
\draw (0,0) -- (1,1.7);
\draw (0,0) -- (2,1.7);
\draw (1,0) -- (0,1.7);
\draw (1,0) -- (2,1.7);
\draw (2,0) -- (0,1.7);
\draw (2,0) -- (1,1.7);
\draw [fill=black] (0,0) circle (1.5pt); 
\draw [fill=black] (1,0) circle (1.5pt); 
\draw [fill=black] (2,0) circle (1.5pt); 
\draw [fill=black] (0.5,.85) circle (1.5pt); 
\draw [fill=black] (1,.85) circle (1.5pt); 
\draw [fill=black] (1.5,.85) circle (1.5pt); 
\draw [fill=black] (0,1.7) circle (1.5pt); 
\draw [fill=black] (1,1.7) circle (1.5pt); 
\draw [fill=black] (2,1.7) circle (1.5pt); 
\end{tikzpicture}
\end{center}
Among the four matroids pictured above, where the bases are given by all triples of points not on a line, the first is algebraic over $k$ if and only if the characteristic of $k$ is $2$, the second and the third are algebraic over any field, and the fourth is algebraic over $k$ if and only if $k$ has nonzero characteristic \cite[Appendix]{Oxley}.
\end{example}

Let $\mathbb{H}^d_n$ be the vector space of all homogeneous polynomials of degree $d$ in $n$ variables with real coefficients, and set
\[
\mathbb{L}^1_n =\{\text{linear forms in $n$ variables with nonnegative coefficients}\}.
\]
Let  $\mathbb{L}^2_n \subseteq \mathbb{H}^2_n$ be the closed subset of quadratic forms with nonnegative coefficients whose Hessians have at most one positive eigenvalue. 
For $d$ larger than $2$, we define $\mathbb{L}^d_n \subseteq \mathbb{H}^d_n$ by setting
\[
\mathbb{L}^d_n=\Big\{f \in \mathbb{M}^d_n\mid \text{$\partial_i f \in \mathbb{L}^{d-1}_n$ for all $i=1,\ldots,n$} \Big\},
\]
where $\mathbb{M}^d_n \subseteq \mathbb{H}^d_n$ is the set of polynomials with nonnegative coefficients whose supports are $\mathrm{M}$-convex.\footnote{By definition, the \emph{support} of a polynomial $f$ is the set of all monomials appearing in $f$ with nonzero coefficients.}
The following characterization in \cite[Theorem 2.25]{BrandenHuh} is central to the theory of Lorentzian polynomials.

\begin{theorem}\label{MainTheorem}
$\mathbb{L}^d_n$  is the set of Lorentzian polynomials in $\mathbb{H}^d_n$.
\end{theorem}

Theorem \ref{MainTheorem} can be used to show that a given polynomial is not Lorentzian. For example, the following polynomials are not Lorentzian because their supports are not $\mathrm{M}$-convex:
\[
x_1^3+x_2^3,  \quad x_1^2x_3+x_2^3.
\]
Note that, in each case, all the partial derivatives $\partial_i f$ are Lorentzian.

One can also use Theorem \ref{MainTheorem} to show that a given polynomial is Lorentzian.
For example, the elementary symmetric polynomial of degree $d$ in $n$ variables is Lorentzian because its support is $\mathrm{M}$-convex and all its quadratic partial derivatives have the Hessian
\[
\left(\begin{array}{cccccc}
0&1&1& \cdots &1 \\
1&0&1& \cdots & 1 \\
1&1&0 & \cdots &1\\
\vdots &\vdots&\vdots&\ddots&\vdots\\
1&1&1&\cdots&0\\
\end{array}
\right),
\]
which has exactly one positive eigenvalue $n-d+1$.
When $n=2$, Theorem \ref{MainTheorem} specializes to the explicit description of bivariate Lorentzian polynomials given in Example~\ref{ExampleBivariate}.

\section{Volume polynomials in projective geometry}\label{SectionProjective}

\subsection{}
The analogous volume polynomial in algebraic geometry is defined as follows: Let $D=(D_1,\ldots,D_n)$ be a collection of semiample divisors on a $d$-dimensional projective variety $Y$ over a field $k$.\footnote{Throughout this paper, a variety over $k$ is by definition a reduced and irreducible scheme of finite type over $k$. A Cartier divisor on a complete variety is \emph{semiample} if some positive multiple moves in a basepoint-free linear system. For background and any undefined terms concerning divisors on varieties and their intersections, we refer to \cite{Lazarsfeld1} and \cite{Fulton}. 
For any field k, the notions of volume polynomial and covolume polynomial over $k$ coincide with those over its algebraic closure $\overline{k}$ \cite[Remark 1.14]{GHMSW}.} The \emph{volume polynomial} of $D$ is 
\[
f_D(x)\coloneqq \frac{1}{d!}\int_Y \bigg(\sum_{i=1}^n x_i D_i\bigg)^d,
\]
which is a homogeneous polynomial of degree $d$ in $x=(x_1,\ldots,x_n)$. 
When the base is the field of complex numbers and each $D_i$ is ample, the restriction of $f_D$ to the positive orthant measures the volume of $Y$ with respect to the K\"ahler class determined by $x$. 

\begin{definition}\label{def:volpoly}
Let $k$ be a field.
\begin{enumerate}[(1)]\itemsep 5pt
\item A homogeneous polynomial $f$ is a \emph{realizable volume polynomial over $k$} if $f=\lambda f_D$ for some $\lambda \in \QQ_{\ge 0}$ and a collection of semiample divisors $D$ on a projective variety $Y$ over $k$. 
\item A homogeneous polynomial $f$ is a \emph{volume polynomial over $k$} if it is a limit of realizable volume polynomials over $k$.
\end{enumerate}
We write $\mathbb{V}^d_n(\mathbb{Q},k)$ for the set of realizable volume polynomials over $k$ of degree $d$ in $n$ variables, and $\mathbb{V}^d_n(\RR,k)$ for the set of all volume polynomials over $k$ of degree $d$ in $n$ variables.
\end{definition}

Recall that a quadratic form 
is Lorentzian if and only if all its coefficients are nonnegative and its Hessian has at most one positive eigenvalue.
It is easy to see that, for any $n$ and any $k$, we have
\begin{align*}
\mathbb{V}^1_n(\QQ,k)&=\big\{\emph{linear forms in $n$ variables with nonnegative rational coefficients}\big\},\\
\mathbb{V}^1_n(\RR,k)&=\big\{\emph{linear forms in $n$ variables with nonnegative coefficients}\big\}.
\end{align*}
By \cite[Theorem 1.8]{HHMWW}, for any $n$ and any $k$, we have
\begin{align*}
\mathbb{V}^2_n(\QQ,k)&=\big\{\emph{Lorentzian quadratic forms in $n$ variables with rational coefficients}\big\},\\
\mathbb{V}^2_n(\RR,k)&=\big\{\emph{Lorentzian quadratic forms in $n$ variables}\big\}.
\end{align*}
Also, by \cite[Theorem 21]{HuhChromatic}, for any $d$ and any $k$, we have
\begin{align*}
\mathbb{V}^d_2(\QQ,k)&=\big\{\emph{Lorentzian bivariate forms of degree $d$ with rational coefficients}\big\},\\
\mathbb{V}^d_2(\RR,k)&=\big\{\emph{Lorentzian bivariate forms of degree $d$}\Big\}.
\end{align*}
In general, $\mathbb{V}^d_n(\RR,k)$ is preserved under a nonnegative linear change of coordinates, and
the normalized coefficients $p_\alpha$ of its members satisfy the \emph{Khovanskii--Teissier inequality}:
\[
p_{\alpha+\mathrm{e}_i-\mathrm{e}_j}p_{\alpha-\mathrm{e}_i+\mathrm{e}_j } \le p_\alpha^2 \ \ \text{for any $\alpha \in \Delta^d_n$ and any $1 \le i<j \le n$.}
\]
As observed in \cite[Section 4.2]{BrandenHuh}, it follows that
\[
\big\{\text{volume polynomials over $k$}\big\} \subseteq \big\{\text{Lorentzian polynomials}\big\} \ \text{for any $k$}.
\]
The \emph{realization problem} for volume polynomials over $k$ is to  identify the volume polynomials over $k$ among Lorentzian polynomials.
The distinction between the notions of realizable volume polynomials and volume polynomials will be relevant in applications to algebraic matroids in Section~\ref{SectionPolymatroid}.

\begin{example}
A standard construction in toric geometry shows that any volume polynomial of rational convex polytopes arises as the volume polynomial of semiample divisors on projective varieties \cite[Section 5.4]{FultonToric}. 
Thus, we have
\[
\big\{\text{volume polynomials of $n$ rational polytopes in $\mathbb{R}^d$}\big\} \subseteq \mathbb{V}^d_n(\mathbb{Q},k) \  \text{for any $k$}.
\]
Since any convex body is a limit of a sequence of rational convex polytopes \cite[Section 1.8]{Schneider}, we have
\[
\big\{\text{volume polynomials of $n$ convex bodies in $\mathbb{R}^d$}\big\} \subseteq \mathbb{V}^d_n(\mathbb{R},k) \  \text{for any $k$}.
\]
As noted in Example~\ref{ExampleQuadratic}, the elementary symmetric polynomial $x_1x_2+x_1x_3+x_1x_4+x_2x_3+x_2x_4+x_3x_4$ is not in the left-hand side.
On the other hand, any Lorentzian quadratic form  is a volume polynomial over $k$ for any $k$, so the inclusion is strict.
\end{example}

\begin{example}
 Let $P=(P_1,\ldots,P_n)$ be a collection of $d \times d$ positive semidefinite Hermitian matrices.
The \emph{volume polynomial} of $P$ is 
\[
f_P(x)\coloneqq \det (x_1P_1+\cdots+x_nP_n),
\]
which is a homogeneous polynomial of degree $d$ in $x=(x_1,\ldots,x_n)$. Using the abelian variety $\mathbb{C}^d/(\mathbb{Z}^d+\mathbb{Z}^d \sqrt{-1})$, one can  show that
\[
\ \big\{\text{volume polynomials of $n$ positive semidefinite $d \times d$ matrices}\big\} \subseteq \mathbb{V}^d_n(\mathbb{R},\mathbb{C}).
\]
Since volume polynomials of positive semidefinite Hermitian matrices are \emph{stable} \cite[Proposition 2.1]{Wagner}, the bivariate Lorentzian polynomial 
$x_1^3+6x_1^2x_2+6x_1x_2^2+2x_2^3$ is not in the left-hand side.
On the other hand, any Lorentzian bivariate form  is a volume polynomial over $k$ for any $k$, so the inclusion is strict. For a detailed discussion of volume polynomials of positive semidefinite Hermitian matrices, see \cite{HMWX}.
\end{example}

\begin{example}
The \emph{Fano matroid} $F_7$ is the rank $3$ matroid on $7$ elements whose bases are all the triples that are not colinear in the following configuration:
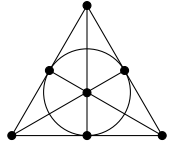
\begin{center}
\begin{tikzpicture}
\draw (0,0) -- (2,0);
\draw (2,0) -- (1,1.73);
\draw (1,1.73) -- (0,0);
\draw (0,0) -- (1.5,.866);
\draw (2,0) -- (.5,.866);
\draw (1,1.73) -- (1,0);
\draw (1,0.577) circle [radius=0.577];
\draw [fill=black] (0,0) circle (1.5pt); 
\draw [fill=black] (2,0) circle (1.5pt); 
\draw [fill=black] (1,1.73) circle (1.5pt); 
\draw [fill=black] (1.5,.866) circle (1.5pt); 
\draw [fill=black] (.5,.866) circle (1.5pt); 
\draw [fill=black] (1,0) circle (1.5pt); 
\draw [fill=black] (1,0.57) circle (1.5pt); 
\end{tikzpicture}
\end{center}
Its \emph{basis generating polynomial} is the cubic polynomial in seven variables
\[
b_{F_7}(x_1,x_2,x_3,x_4,x_5,x_6,x_7)=\sum_{ijk \in F_7} x_i x_j x_k,
\]
where the sum is over all the $28$ bases of the Fano matroid. According to  \cite[Example 5.5]{GHMSW}, we have
\[
b_{F_7} \in \mathbb{V}^3_7(\mathbb{Q},k) \ \ \text{if and only if} \ \ \text{char}(k)=2.
\]
Is $b_{F_7}$ a volume polynomial over $k$ when the characteristic of $k$ is not $2$?
See Conjecture~\ref{ConjectureIndependence} below.
\end{example}

\begin{example}\label{ExampleLinear}
If a matroid $J$ of rank $d$ on $n$ elements is linear over $k$, then its basis generating polynomial $f_J$ is a homogeneous polynomial of degree $d$ in $n$ variables.
The \emph{arrangement Schubert variety} of any linear realization of $J$ over $k$ witnesses  the fact that $f_J$ is a realizable volume polynomial over $k$ \cite[Section 1.3]{BHMPW}.  
\end{example}

\begin{example}\label{ExampleSchur}
The \emph{Schur module} $\mathrm{V}(\lambda)$ of a Young diagram $\lambda$ is  the irreducible representation of  the general linear group  $\mathrm{GL}_n(\mathbb{C})$ with highest weight $\lambda$.
It  has the weight space decomposition
\[
\mathrm{V}(\lambda)=\bigoplus_\mu \mathrm{V}(\lambda)_\mu \ \ \text{with} \ \  \dim \mathrm{V}(\lambda)_\mu=K_{\lambda\mu},
\]
where $K_{\lambda\mu}$ is the \emph{Kostka number} counting semistandard Young tableaux of given shape $\lambda$ and  weight $\mu$ \cite[Section 8.3]{FultonYoung}.
The \emph{normalized Schur polynomial} is the generating polynomial
\[
f_\lambda(x_1,\ldots,x_n) \coloneqq \sum_{\mu} K_{\lambda\mu}\hspace{0.3mm} x^{[\mu]}, \ \ \text{where} \ \  x^{[\mu]} \coloneqq \frac{x^\mu}{\mu!}= \frac{x_1^{\mu_1}}{\mu_1!} \cdots  \frac{x_n^{\mu_n}}{\mu_n!}.
\]
For example, for the Young diagram $\lambda=\raisebox{-1.0ex}{\scalebox{0.7}{\yng(2,1)}}$, we have
\[
f_\lambda(x_1,x_2,x_3)=\frac{1}{2}x_1^2x_2 +\frac{1}{2}x_1^2x_3+\frac{1}{2}x_1x_2^2+\frac{1}{2}x_2^2 x_3+\frac{1}{2}x_1 x_3^2+\frac{1}{2}x_2 x_3^2 +2x_1x_2x_3.
\]
The proof of \cite[Theorem 3]{HMMS} shows that any normalized Schur polynomial is a  realizable volume polynomial over $k$ for any $k$.
\end{example}

\begin{example}
The product of two realizable volume polynomials over $k$ is a realizable volume polynomial over $k$: If $f_1(x)$ is the realizable volume polynomial obtained from a collection of semiample divisors $D_1$ on $Y_1$
and $f_2(x)$ is the realizable volume polynomial obtained from a collection of semiample divisors $D_2$ on $Y_2$,
then $f_1(x)f_2(x)$ is the realizable volume polynomial obtained from the collection of  semiample divisors $\pi_1^*D_1+\pi_2^*D_2$ on $Y_1 \times Y_2$, where the addition is defined componentwise.
\end{example}

It is known that $\mathbb{V}^d_n(\QQ,k)$, and hence $\mathbb{V}^d_n(\RR,k)$, depends only on the characteristic of $k$ \cite[Proposition 2.10]{GHMSW}.
It is not known whether $\mathbb{V}^d_n(\RR,k)$ depends on $k$
when $d \ge 3$ and $n \ge 3$. 

\begin{conjecture}\label{ConjectureIndependence}
The set of volume polynomials over $k$ is independent of the choice of $k$.
\end{conjecture}

\subsection{}
As observed in  \cite[Example~14]{HuhICM}, we have the proper inclusion
\[
\mathbb{V}^d_n(\mathbb{R},k)\subsetneq \mathbb{L}^d_n \ \ \text{for any field $k$ when $d \ge 3$ and $n \ge 3$.}
\]
Thus, the realization problem for volume polynomials over $k$ has a nontrivial answer in these cases.
For example, consider the cubic polynomial
\[
f=14x_1^3+6x_1^2x_2+24x_1^2x_3+12x_1x_2x_3+6x_1x_3^2+3x_2x_3^2.
\]
The support of $f$ is $\mathrm{M}$-convex, as it is the set of all lattice points of the following integral generalized permutohedron:
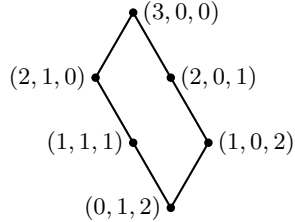
\[
\begin{tikzpicture}[scale=1.00, every node/.style={font=\small}]
  \def\h{0.8660254038} 

  \newcommand{\pt}[3]{({(#2-#3)/2},{-\h*(#2+#3)})}

  \draw[thick]
    \pt{3}{0}{0} --
    \pt{2}{0}{1} --
    \pt{0}{2}{1} --
    \pt{1}{2}{0} -- cycle;

  \fill \pt{3}{0}{0} circle (1.6pt);
  \node[right]      at \pt{3}{0}{0} {$(3,0,0)$};

  \fill \pt{2}{1}{0} circle (1.6pt);
  \node[right]       at \pt{2}{1}{0} {$(2,0,1)$};

  \fill \pt{2}{0}{1} circle (1.6pt);
  \node[left]      at \pt{2}{0}{1} {$(2,1,0)$};

  \fill \pt{1}{2}{0} circle (1.6pt);
  \node[right] at \pt{1}{2}{0} {$(1,0,2)$};

  \fill \pt{0}{2}{1} circle (1.6pt);
  \node[left] at \pt{0}{2}{1} {$(0,1,2)$};

  \fill \pt{1}{1}{1} circle (1.6pt);
  \node[left]      at \pt{1}{1}{1} {$(1,1,1)$};

\end{tikzpicture}
\]
The Hessians of the partial derivatives
 $\partial_1 f, \partial_2 f, \partial_3 f$ are
 \[
 \begin{pmatrix}
84 & 12 & 48 \\
12 & 0 & 12 \\
48 & 12 & 12
\end{pmatrix}, \ \ 
\begin{pmatrix}
12 & 0 & 12 \\
0 & 0 & 0 \\
12 & 0 & 6
\end{pmatrix}, \ \ 
\begin{pmatrix}
48 & 12 & 12 \\
12 & 0 & 6 \\
12 & 6 & 0
\end{pmatrix}, 
 \]
 each of which has exactly one positive eigenvalue.
 Then, by Theorem~\ref{MainTheorem},   $f$ is a Lorentzian polynomial.
The fact that $f$ is not a volume polynomial over $k$ follows from the \emph{reverse Khovanskii--Teissier inequality} \cite[Theorem 5.7]{LehmannXiao}:
For any nef divisors $D_1,D_2,D_3$ on a $d$-dimensional  projective variety $Y$ and any $e \le d$, 
\[
\binom{d}{e} \left( \int_Y D_1^{d-e} D_2^e\right)  \left( \int_Y  D_1^e  D_3^{d-e}  \right) \ge \left( \int_Y  D_1^d \right) \left( \int_Y  D_2^e  D_3^{d-e}\right).
\]
The complex analytic proof of the inequality in \cite{LehmannXiao} relies on the Calabi--Yau theorem \cite{Yau}. The algebraic proof of the inequality in \cite{JiangLi} using Okounkov bodies works over any algebraically closed field. 
Since the notion of volume polynomial over $k$ coincides with that over its algebraic closure $\overline{k}$ \cite[Remark 1.14]{GHMSW}, the inequality remains valid for any $k$. 
As mentioned before, 
\[
\big\{\text{volume polynomials of $n$ convex bodies in $\mathbb{R}^d$}\big\} \subseteq \mathbb{V}^d_n(\mathbb{R},k) \  \text{for any $k$}.
\]
Since Lorentzian polynomials are strongly log-concave \cite[Theorem 2.31]{BrandenHuh}, the Lorentzian cubic $f$ provides a counterexample to Gurvits' conjecture that a strongly log-concave homogeneous polynomial in three variables with nonnegative coefficients is the volume polynomial of three convex bodies
\cite[Conjecture 4.1]{GurvitsL}.

In \cite{HMWX}, the authors introduce a new family of inequalities for volume polynomials that subsumes both the Khovanskii--Teissier and the reverse Khovanskii--Teissier inequalities as special cases.

\subsection{}

The set of volume polynomials over $k$  is preserved under any nonnegative linear change of coordinates, as in the case of volume polynomials of convex bodies (Section~\ref{SectionConvex}).
More precisely, for any  $n \times m$ matrix $A$ with nonnegative rational entries and sets of variables $x=(x_1,\ldots,x_n)$ and $y=(y_1,\ldots,y_m)$, we have the implication
\[
f(x) \in \mathbb{V}^d_n(\mathbb{Q},k) \ \ \Longrightarrow \ \ f(Ay) \in \mathbb{V}^d_m(\mathbb{Q},k).
\]
It follows that, for any  $n \times m$ matrix $A$ with nonnegative real entries, we have
\[
f(x) \in \mathbb{V}^d_n(\mathbb{R},k) \ \ \Longrightarrow \ \ f(Ay) \in \mathbb{V}^d_m(\mathbb{R},k).
\]
Here are two additional basic operations on the set of realizable volume polynomials over $k$, and hence on the set of volume polynomials over $k$.
Let $E$ be a finite set indexing the variables in a polynomial ring with real coefficients. 

\begin{definition}\label{DefinitionMinor}
For a nonzero degree $d$ homogeneous polynomial  $f \in \mathbb{R}[x_i]_{i \in E}$ and an element $j \in E$, we write
\[
f=\sum_{e=e_{\min}}^{e_{\max}} f_{d-e} \frac{x_j^e}{e!}, 
\]
where $f_{d-e}$ are polynomials in $ \mathbb{R}[x_i]_{i \neq j}$ that are nonzero for $e=e_{\min}, e_{\max}$.
\begin{enumerate}[(1)]\itemsep 5pt
\item The \emph{deletion} of $f$ by $j$ is the degree $d$ homogeneous polynomial 
\[
f\setminus j=\sum_{e=e_{\min}}^{e_{\max}-1} f_{d-e} \frac{x_j^e}{e!}.
\]
\item The \emph{contraction} of $f$ by $j$ is the degree $d-1$ homogeneous polynomial
\[
f/j=\sum_{e=e_{\min}+1}^{e_{\max}} f_{d-e} \frac{x_j^{e-1}}{(e-1)!}.
\]
\end{enumerate}
A \emph{minor} of $f$ is a polynomial obtained from $f$ by a sequence of deletion and contraction operations.
\end{definition}

When applied to the spanning tree polynomials of graphs and the basis generating polynomials of matroids, Definition~\ref{DefinitionMinor} recovers the corresponding notions of contraction, deletion, and minor in the context of graph theory and matroid theory \cite[Chapter 3]{Oxley}.\footnote{In \cite[Section 1.2]{Oxley}, a matroid is required to have at least one basis.  This leads to minor differences in the definition of $f \setminus j$ and $f / j$ when $j$ is a loop (when $j$ is contained in no basis of the matroid) or a coloop (when $j$ is contained in every basis of the matroid).}
For a discussion of minors in the more general framework of polymatroids over tracts, see \cite[Section 2.2]{BHKL-polymatroids}.

A special case of  \cite[Corollary 3.3]{{GHMSW}} implies that  deletions and contractions of realizable volume polynomials over $k$ are realizable volume polynomials over $k$. 
For a systematic study of linear operators preserving the set of realizable volume polynomials over $k$, see Section~\ref{SectionVolumeCovolume}.

\begin{proposition}\label{PropMinor}
 Any minor of a  realizable volume polynomial over $k$ is a  realizable volume polynomial over $k$.
\end{proposition}

It follows that any minor of a volume polynomial over $k$ is a volume polynomial over $k$. 
This contrasts  with the fact that the set of volume polynomials of convex bodies is not closed under minors.

\begin{example}\label{ExampleMinor}
Let  $C_1,C_2,C_3,C_4$ be four equiangular unit segments in $\mathbb{R}^3$, and let $C_5$ be the unit ball in $\mathbb{R}^3$. 
The volume polynomial for $C=(C_1,C_2,C_3,C_4,C_5)$ is the cubic in five variables
\[
f_C=\frac{4 \pi}{3} x_5^3+\pi  \Bigg[\sum_{1 \le  i \le 4} x_i \Bigg] x_5^2+\frac{4\sqrt{2}}{3}  \Bigg[\sum_{1 \le  i<j \le 4} x_ix_j \Bigg] x_5  +\frac{4 \sqrt{3}}{9} \Bigg[\sum_{1 \le  i<j<k \le 4} x_ix_jx_k \Bigg].
\]
Recall from Example~\ref{ExampleQuadratic}  that
the quadratic elementary symmetric polynomial in $x_1,x_2,x_3,x_4$ is not a volume polynomial of convex bodies.
However, it is a minor of the volume polynomial $f_C$.
\end{example}

To what extent do the volume polynomials arising in algebraic geometry coincide with those arising in convex geometry?
Proposition~\ref{PropMinor} shows that any minor of a volume polynomial of convex bodies is a volume polynomial over $k$ for any $k$.
The following strengthening of Conjecture~\ref{ConjectureIndependence} was suggested during a discussion with Shouda Wang.

\begin{conjecture}
Every volume polynomial over $k$ is a limit of minors of volume polynomials of convex bodies. 
\end{conjecture}

\subsection{}\label{SectionAnalytic}

A real $(1,1)$-class $[\omega]$ on a compact K\"ahler manifold $Y$ is \emph{semipositive}
 if it contains a smooth semipositive representative, that is, if there is a smooth function $\varphi$ on $Y$ such that
\[
\omega+i\partial \overline{\partial} \varphi \ge 0.
\]

\begin{definition}\label{DefinitionAnalytic}
A degree $d$ homogeneous polynomial $f$ in $n$ variables  is a \emph{realizable analytic volume polynomial} if there is a $d$-dimensional compact K\"ahler manifold $Y$ and semipositive classes $[\omega_1],\dots,[\omega_n]$ such that 
     \[
     f(x_1,\dots,x_n)=\frac{1}{d!}\int_Y (x_1\omega_1+\dots +x_n\omega_n)^{\wedge d}.
     \]
    A homogeneous polynomial $f$ is an \emph{analytic volume polynomial} if it is a limit of realizable analytic volume polynomials.
    \end{definition}
    
By \cite{Gromov}, the $(1,1)$-part of the cohomology of $Y$ satisfies the \emph{mixed Hodge--Riemann relations}, and hence
\[
 \big\{\text{analytic volume polynomials}\big\} \subseteq \big\{\text{Lorentzian polynomials}\big\}.
\]
It follows that the support of an analytic volume polynomial is $\mathrm{M}$-convex, defining the class of \emph{analytic polymatroids} \cite[Section 5]{GHMSW}. 
On the other hand, the resolution of singularities for complex projective varieties implies that any realizable volume polynomial over $\mathbb{C}$ is a 
realizable analytic volume polynomial, and hence
\[
\big\{\text{volume polynomials over $\mathbb{C}$}\big\} \subseteq \big\{\text{analytic volume polynomials}\big\}.
\]
The answers to the following basic questions regarding analytic volume polynomials remain unknown.

\begin{question}
Is there an analytic volume polynomial that is not a volume polynomial over $\mathbb{C}$?
\end{question}

\begin{question}
Is the class of analytic volume polynomials closed under taking minors?
\end{question}








\section{Linear operators preserving volume polynomials}\label{SectionVolumeCovolume}

\subsection{}

The set of volume polynomials over $k$ is, in a precise sense, dual to the set of \emph{covolume polynomials} over $k$. 
To define covolume polynomials and state their main properties, it will be convenient to work with the dual pair of polynomial rings
\[
\RR[\partial]= \RR[\partial_i]_{i \in E}  \quad \text{and} \quad 
\RR[x]= \RR[x_i]_{i \in E}.
\]
We write $ \mathbb{Z}^E_{\ge 0}$ for the set of exponent vectors of the monomials in the two polynomial rings, and set
\[
\partial^\alpha \coloneqq \prod_{i \in E}\partial_i^{\alpha_i}  \ \ \text{and} \ \  x^{[\alpha]} \coloneqq \prod_{i \in E}\frac{x_i^{\alpha_i}}{\alpha_i!}  \ \ \text{for $\alpha \in \mathbb{Z}^E_{\ge 0}$.}
\]
The polynomial ring $\RR[\partial]$ acts on $\RR[x]$ as differential operators by the usual rule
\[
\partial^\alpha \circ x^{[\beta]}\coloneqq \begin{cases} x^{[\beta-\alpha]}& \text{if $\alpha \le \beta$,}\\ \hfill 0\hfill & \text{if otherwise,}\end{cases}
\]
where $\alpha \le \beta$ means that their components satisfy $\alpha_i \le \beta_i$ for  all $i \in E$.
For any further conventions for multivariate polynomials, we refer to \cite[Section~2]{BrandenHuh}.
For $\mu \in \ZZ^E_{\ge 0}$, we consider 
\[
\RR[\partial]_{\le \mu}\coloneqq \textrm{span}(\partial^\alpha)_{\alpha \le \mu} \quad \text{and} \quad
\RR[x]_{\le \mu}\coloneqq \textrm{span}(x^\alpha)_{\alpha \le \mu}.
\]
Then $\RR[x]_{\le \mu}$ is an $\RR[\partial]$-submodule of $\RR[x]$ generated by $x^{[\mu]}$, and the linear map 
\[
\RR[\partial]_{\le \mu} \longrightarrow \RR[x]_{\le \mu}, \qquad \partial^\alpha \longmapsto \partial^\alpha \circ x^{[\mu]}=x^{[\mu-\alpha]}
\]
is an isomorphism of finite-dimensional vector spaces.

\begin{definition}\label{DefinitionCovolume}
Let $g$ be a homogeneous polynomial in $\RR[\partial]_{\le \mu}$.
\begin{enumerate}[(1)]\itemsep 5pt
\item We say that $g$ is a \emph{realizable covolume polynomial over $k$} if $g(\partial) \circ x^{[\mu]}$ is a realizable volume polynomial over $k$.
\item We say that  $g$ is a \emph{covolume polynomial over $k$} if it is a limit of realizable covolume polynomials over $k$.
\end{enumerate}
\end{definition}

As observed in \cite[Remark 2.2]{Aluffi}, the property of being a realizable covolume polynomial over $k$ does not depend on the choice of $\mu$. 
This follows from the translation invariance
\begin{multline*}
\Big(\text{ $\sum_\alpha c_\alpha x^{[\alpha]}$ is a realizable volume polynomial over $k$}\Big) \Longleftrightarrow \\
\Big(\text{ $\sum_\alpha c_\alpha x^{[\alpha+\beta]}$ is a realizable volume polynomial over $k$}\Big), \ \text{for any $\beta \in \mathbb{Z}^E_{\ge 0}$.}
\end{multline*}
For the cone construction that justifies this, see  \cite[Section 2]{GHMSW}.

\begin{remark}
According to \cite[Section 2.1]{BHKL-polymatroids}, 
the \emph{dual} of an $\mathrm{M}$-convex set $J \subseteq \mathbb{Z}^E_{\ge 0}$, defined up to translation in $\mathbb{Z}^E_{\ge 0}$, is the  $\mathrm{M}$-convex subset 
\[
\mu-J  \coloneqq \{\mu-\alpha \mid \alpha \in J \}\subseteq \mathbb{Z}^E_{\ge 0},
\] 
where $\mu$ is any nonnegative integral vector satisfying $\alpha \le \mu$ for all $\alpha \in J$.\footnote{A standard choice used in \cite[Section 2.1]{BHKL-polymatroids} is to take $\mu$ to be the \emph{duality vector} $\delta_J=\delta^+_J+\delta^-_J$, where $\delta^+_J=\sup J$ and $\delta^-_J=\inf J$ under the partial order on $\mathbb{Z}^E_{\ge 0}$ defined by $\alpha \le \beta$ if and only if $\alpha_i \le \beta_i$ for all $i \in E$.}
Since 
the support of a covolume polynomial over $k$ is the dual of the support of a volume polynomial over $k$,
the support of a covolume polynomial is an $\mathrm{M}$-convex set.
\end{remark}

\begin{remark}\label{RemarkMinors}
The class of (realizable) volume polynomials over $k$ is closed under nonnegative (rational) linear changes of coordinates, as well as under taking products and minors.
Similarly, the class of (realizable) covolume polynomials over $k$ is closed under nonnegative (rational) linear changes of coordinates \cite[Theorem 2.7]{GHMSW}, as well as under taking products and minors \cite[Theorem 1.5]{GHMSW}.
It is interesting to note that the corresponding statements for volume polynomials and covolume polynomials sometimes have substantially different proofs.
\end{remark}

\begin{example}
Let $f_J$ be  the basis generating polynomial of a matroid $J$ linear over $k$.
By Example~\ref{ExampleLinear}, $f_J$ is a realizable volume polynomial over $k$.  
 Since the dual of a linear matroid over $k$ is linear over the same field \cite[Corollary 2.2.9]{Oxley},  $f_J$ is a realizable covolume polynomial over $k$ as well.  
\end{example}

\begin{example}
The proof of \cite[Theorem 6]{HMMS}  shows that the \emph{Schubert polynomial} $s_w(\partial)$ is a realizable covolume polynomial over $k$ for any permutation $w$ and any $k$.
In particular, any  \emph{Schur polynomial} is a realizable covolume polynomial over $k$ for any field $k$.
By Example~\ref{ExampleSchur}, any \emph{normalized Schur polynomial} is a realizable volume polynomial over $k$ for any $k$. 
It is not known whether \emph{normalized Schubert polynomials} are realizable volume polynomials over $k$ for any $k$. See \cite[Conjecture 15]{HMMS} for a weaker statement.
\end{example}

\begin{example}
The convex polytope in Example~\ref{ExampleVolumeCovolume} shows that, for any $k$, 
\begin{multline*}
x_1x_2x_3+x_1x_2x_4+x_1x_2x_5+x_1x_3x_4+x_1x_3x_5 \\
+x_1x_4x_5+x_2x_3x_4+x_2x_3x_5+x_2x_4x_5 +4x_3x_4x_5
\end{multline*}
is a realizable volume polynomial over $k$. It follows that, for any $k$, 
\[
\partial_4\partial_5+\partial_3\partial_5+\partial_3\partial_4+\partial_2\partial_5+\partial_2\partial_4+\partial_2\partial_3+ \partial_1\partial_5+ \partial_1 \partial_4+ \partial_1 \partial_3+4 \partial_1\partial_2
\]
is a realizable covolume polynomial over $k$.
This covolume polynomial over $k$ is not a Lorentzian polynomial, and hence it is not a volume polynomial over $k$ for any $k$.
\end{example}

Conjecture~\ref{ConjectureCovolume} suggests that a quadratic polynomial  $\sum_{1 \le i<j \le n} q_{ij} \partial_i \partial_j$ is a covolume polynomial over $k$ if and only if
its coefficients are nonnegative and satisfy the \emph{Pl\"ucker relations} over the triangular hyperfield $\mathbb{T}_2$:
\[
\sqrt{q_{ij}q_{kl}} \le  \sqrt{q_{ik}q_{jl}} + \sqrt{q_{il}q_{jk}} \ \ \text{for any $1 \le i<j<k<l\le n$.}
\]
For general discussions of Grassmannians over triangular hyperfields, see \cite{BHKL-triangular,BHKL-triangular2}. 

\subsection{}

The main result of \cite{GHMSW} is the following characterization of realizable covolume polynomials. This parallels the characterization of dually Lorentzian polynomials in \cite[Theorem 1.2]{RSW}.

\begin{theorem}\label{thm:covolumecharacterization}
The following conditions are equivalent for any $g \in \QQ[\partial]$.
\begin{enumerate}[(1)]\itemsep 5pt
\item The polynomial $g$ is a realizable covolume polynomial over $k$.
\item For any realizable volume polynomial $f$ over $k$, the polynomial $g(\partial) \circ f(x)$ is a realizable volume polynomial over $k$.
\end{enumerate}
\end{theorem}

The corresponding characterization of covolume polynomials is  obtained by taking limits.

\begin{corollary}\label{cor:covolumecharacterization}
The following conditions are equivalent for any $g \in \RR[\partial]$.
\begin{enumerate}[(1)]\itemsep 5pt
\item The polynomial $g$ is a covolume polynomial over $k$.
\item For any volume polynomial $f$ over $k$, the polynomial $g(\partial) \circ f(x)$ is a volume polynomial over $k$.
\end{enumerate}
\end{corollary}

Example~\ref{ExampleMinor} shows that the set of volume polynomials of convex bodies does not satisfy the analogous statement.\footnote{For example, $\partial_5 f_C$ is not a volume polynomial of convex bodies.}
For a parallel statement characterizing volume polynomials as linear operators preserving covolume polynomials, see \cite[Theorem 1.9]{GHMSW}. 

Corollary~\ref{cor:covolumecharacterization} can be used to deduce new inequalities for mixed volumes of convex bodies, or more generally, for intersection numbers of nef divisors on a projective variety.
For instance, given a Schubert polynomial $s_w(\partial)$ and a volume polynomial $f_C(x)$, any known inequality for the coefficients of a volume polynomial can be applied to 
$s_w(\partial) \circ f_C(x)$ 
to produce another inequality for the coefficients of $f_C(x)$.
For an overview of known inequalities for the coefficients of the volume polynomial, such as the Khovanskii--Teissier inequality or the reverse Khovanskii--Teissier inequality, see \cite{HMWX}. 

\begin{question}\label{RemarkAnalytic}
One can define \emph{(realizable) analytic covolume polynomials} as the duals of (realizable) analytic volume polynomials as in Definition~\ref{DefinitionCovolume}. Do they satisfy the analogues of Theorem~\ref{thm:covolumecharacterization} and Corollary~\ref{cor:covolumecharacterization}?
\end{question}

\subsection{}

The \emph{symbol theorem} for Lorentzian polynomials states that, if the symbol of a linear operator $T$ is a Lorentzian polynomial, then $T$ sends Lorentzian polynomials to Lorentzian polynomials  \cite[Theorem 3.2]{BrandenHuh}. 
Theorem~\ref{thm:covolumecharacterization} can be used to derive a volume polynomial analogue of the symbol theorem for Lorentzian polynomials.\footnote{The study of symbols of linear operators dates back to G\r{a}rding \cite{Garding} and appears prominently in the work of Borcea and Br\"and\'en on the P\'olya--Schur program for stable polynomials \cite{BorceaBranden2,BorceaBranden1}.}

Let $x=(x_i)_{i \in E}$ and $y=(y_j)_{j \in F}$ be two finite sets of variables. 
Let $T$ be a  homogeneous linear operator\footnote{This means that $T$ is linear over $\mathbb{R}$ and $\deg\, T(x^\alpha)-\deg\, x^\alpha \in \ZZ$ does not depend on $\alpha \le \mu$.} 
\[
T: \RR[x]_{\le \mu} \longrightarrow \RR[y]_{\le \nu}, \ \  \text{where $\mu \in \mathbb{Z}^E_{\ge 0}$ and $\nu \in \mathbb{Z}^F_{\ge 0}$}. 
\]
The \emph{symbol} of $T$ is the homogeneous polynomial in variables $(x,y)$ given by
\[
\textrm{sym}_T(x,y)=\sum_{0 \le \alpha \le \mu}   T(x^{[\alpha]}) x^{[\mu-\alpha]}.
\]

\begin{theorem}\label{thm:symbol}
If the symbol of $T$ is a realizable volume polynomial over $k$, then $T$ sends realizable volume polynomials over $k$ to realizable volume polynomials over $k$.
\end{theorem}

\begin{corollary}
If the symbol of $T$ is a volume polynomial over $k$, then $T$ sends volume polynomials over $k$ to volume polynomials over $k$.
\end{corollary}

Geometrically, one may view $T$ as a graded linear map between Chow groups
\[
\varphi_T:\textrm{CH}(\mathbb{P}^\mu) \otimes \mathbb{R} \to \textrm{CH}(\mathbb{P}^\nu) \otimes \mathbb{R}, \ \ \text{where $\mathbb{P}^\mu = \prod_{i \in E} \mathbb{P}^{\mu_i}$ and $\mathbb{P}^\nu = \prod_{j \in F} \mathbb{P}^{\nu_j}$.}
\]
If this map is induced by an irreducible correspondence $\Gamma \subseteq \mathbb{P}^\mu \times \mathbb{P}^\nu$ so that
    \[
    \varphi_{T}(\Lambda)=p_{2*}\big(\Gamma\cap p_1^*(\Lambda)\big), 
    \]
then, by  \cite[Lemma 2.1]{GHMSW}, it preserves the classes of irreducible cycles up to a rational multiple.

The symbol theorem for realizable volume polynomials shows that many familiar operators from the theory of Lorentzian polynomials preserve realizable volume polynomials over $k$ for any $k$:
\begin{enumerate}[(1)]\itemsep 5pt
\item The \emph{upper truncation operators} and the \emph{lower truncation operators} preserve realizable volume polynomials over $k$  \cite[Corollary 3.3]{GHMSW}.
\item The \emph{polarization operator} $\Pi^\uparrow$ preserves realizable volume polynomials over $k$ \cite[Proposition 4.1]{GHMSW}.
\item The \emph{normalization operator} $\mathrm{N}$ preserves realizable volume polynomials over $k$  \cite[Proposition 4.2]{GHMSW}.
\item For any nonnegative rational number $t$, the \emph{interlacing operator} $1+t x_i\partial_j$ preserves realizable volume polynomials over $k$ \cite[Proposition 4.3]{GHMSW}.
\item For any nonnegative rational number $t$, the \emph{symmetric exclusion process} $\Phi^{i,j}_t$ preserves realizable multiaffine volume polynomials over $k$ \cite[Proposition 4.4]{GHMSW}.
\end{enumerate}
The  corresponding statements for volume polynomials over $k$ follow from taking limits.

    
   

\section{Realization problems for polymatroids}\label{SectionPolymatroid}

Recall that the \emph{support} of a polynomial $f$ in $\RR[x_i]_{i \in E}$ is the set of all exponent vectors $\alpha \in \mathbb{Z}^E_{\ge 0}$ such that the monomial $x^\alpha$ appears in $f$ with nonzero coefficient.

\begin{question}
If $f$ is a volume polynomial over $k$ (Definition~\ref{def:volpoly}), then the support of $f$ is the set of bases of a polymatroid.
Which polymatroids arise in this way?
\end{question}

\begin{question}
If $f$ is a covolume polynomial over $k$ (Definition~\ref{DefinitionCovolume}), then the support of $f$ is the set of bases of a polymatroid.
Which polymatroids arise in this way?
\end{question}

\begin{question}
If $f$ is an analytic volume polynomial (Definition~\ref{DefinitionAnalytic}), then the support of $f$ is the set of bases of a polymatroid.
Which polymatroids arise in this way?
\end{question}

\begin{question}
If $f$ is an analytic covolume polynomial (Question~\ref{RemarkAnalytic}), then the support of $f$ is the set of bases of a polymatroid.
Which polymatroids arise in this way?
\end{question}

At present, the author is not aware of any obstructions in any of the above cases.

The following connection between algebraic matroids and the support of \emph{realizable} volume polynomials is known  \cite[Proposition 5.4]{GHMSW}. 
A polymatroid on $E$ is  \emph{algebraic over $k$} if there are field extensions $k \subseteq \ell_i \subseteq \ell$  for $i \in E$
 such that 
 \[
h(A)= \text{trdeg}_k \Big(\vee_{i \in A} \ell_i\Big) \ \ \text{for all $A \subseteq E$},
 \]
 where $h$ is the rank function of the polymatroid.

\begin{proposition}\label{PropositionAlgebraicPolymatroids}
A polymatroid is algebraic over $k$ if and only if it is the support of a realizable volume polynomial over $k$.
\end{proposition}

Proposition~\ref{PropositionAlgebraicPolymatroids} implies that every minor of an algebraic polymatroid over $k$ is an algebraic polymatroid over $k$ \cite[Section 5]{GHMSW}.\footnote{Definition~\ref{DefinitionMinor}, applied to the basis generating polynomials of polymatroids, defines the notion of \emph{minor} of polymatroids  \cite[Section 2.1]{BHKL-polymatroids}.}
In the classical case of matroids, this statement is typically deduced from a theorem of Lindstr\"om \cite{Lindstrom89}, who proved Piff's conjecture 
 that $M$ is algebraic over $k$ if $M$ is algebraic over an extension of $k$; see  \cite[Corollary 6.7.14]{Oxley}. 
Since the set of realizable volume polynomials over $k$ depends only on the characteristic of $k$ \cite[Proposition 2.10]{GHMSW}, Proposition~\ref{PropositionAlgebraicPolymatroids}  gives the following version of Lindstr\"om's theorem for polymatroids.

\begin{corollary}
A polymatroid  is algebraic over some field of characteristic $p$ if and only if 
it is algebraic over all fields of characteristic $p$.
\end{corollary}

Another consequence of Proposition~\ref{PropositionAlgebraicPolymatroids} is that the  intersection $M_1 \wedge M_2$ of algebraic matroids over $k$ is an algebraic matroid over $k$  \cite[Theorem 5.11]{GHMSW}.
  This generalizes Piff's theorem that the truncation of an algebraic matroid is algebraic \cite[Section 11.3]{WelshMatroidTheory}.\footnote{In \cite[Section 11.3]{WelshMatroidTheory}, Welsh proves the dual statement that the union $M_1 \vee M_2$ of algebraic matroids over $k$ is an algebraic matroid over $k$. It is interesting to note that, as in Remark~\ref{RemarkMinors}, the proof of the statement for $M_1 \wedge M_2$ is substantially different from that for $M_1 \vee M_2$.}

Is the support of a realizable covolume polynomial over $k$  an algebraic polymatroid over $k$? 
This question extends the following long-standing open problem in matroid theory. 
For up-to-date discussions, see \cite{BFP,Hochstattler}.

\begin{question}
Is the dual of an algebraic matroid over $k$ algebraic over $k$?
\end{question}


\bibliographystyle{amsalpha}
\bibliography{biblio}

@incollection {Gromov,
    AUTHOR = {Gromov, Mikhael},
     TITLE = {Convex sets and {K}\"ahler manifolds},
 BOOKTITLE = {Advances in differential geometry and topology},
     PAGES = {1--38},
 PUBLISHER = {World Sci. Publ., Teaneck, NJ},
      YEAR = {1990},
      ISBN = {981-02-0494-9; 981-02-0495-7},
   MRCLASS = {52A40 (53C55)},
  MRNUMBER = {1095529},
MRREVIEWER = {Y.\ Tashiro},
}

@article {BorceaBranden1,
    AUTHOR = {Borcea, Julius and Br\"and\'en, Petter},
     TITLE = {P\'olya-{S}chur master theorems for circular domains and their
              boundaries},
   JOURNAL = {Ann. of Math. (2)},
  FJOURNAL = {Annals of Mathematics. Second Series},
    VOLUME = {170},
      YEAR = {2009},
    NUMBER = {1},
     PAGES = {465--492},
      ISSN = {0003-486X,1939-8980},
   MRCLASS = {30C10 (65H05 65J05)},
  MRNUMBER = {2521123},
MRREVIEWER = {Vania\ D.\ Mascioni},
       DOI = {10.4007/annals.2009.170.465},
       URL = {https://doi.org/10.4007/annals.2009.170.465},
}

@article {BorceaBranden2,
    AUTHOR = {Borcea, Julius and Br\"and\'en, Petter},
     TITLE = {The {L}ee-{Y}ang and {P}\'olya-{S}chur programs. {I}. {L}inear
              operators preserving stability},
   JOURNAL = {Invent. Math.},
  FJOURNAL = {Inventiones Mathematicae},
    VOLUME = {177},
      YEAR = {2009},
    NUMBER = {3},
     PAGES = {541--569},
      ISSN = {0020-9910,1432-1297},
   MRCLASS = {47B38 (33C45 82B26)},
  MRNUMBER = {2534100},
MRREVIEWER = {Roland\ K. W. Roeder},
       DOI = {10.1007/s00222-009-0189-3},
       URL = {https://doi.org/10.1007/s00222-009-0189-3},
}

@article{Wagner,
  author  = {David G. Wagner},
  title   = {Multivariate stable polynomials: theory and applications},
  journal = {Bull. Amer. Math. Soc. (N.S.)},
  volume  = {48},
  number  = {1},
  pages   = {53--84},
  year    = {2011},
  doi     = {10.1090/S0273-0979-2010-01321-5},
}

@book {Schneider,
    AUTHOR = {Schneider, Rolf},
     TITLE = {Convex bodies: the {B}runn-{M}inkowski theory},
    SERIES = {Encyclopedia of Mathematics and its Applications},
    VOLUME = {151},
   EDITION = {expanded},
 PUBLISHER = {Cambridge University Press, Cambridge},
      YEAR = {2014},
     PAGES = {xxii+736},
      ISBN = {978-1-107-60101-7},
   MRCLASS = {52-02 (52A20 52A39)},
  MRNUMBER = {3155183},
MRREVIEWER = {Andrea Colesanti},
}

@book {Fulton,
    AUTHOR = {Fulton, William},
     TITLE = {Intersection theory},
    SERIES = {Ergebnisse der Mathematik und ihrer Grenzgebiete. 3. Folge. A
              Series of Modern Surveys in Mathematics [Results in
              Mathematics and Related Areas. 3rd Series. A Series of Modern
              Surveys in Mathematics]},
    VOLUME = {2},
   EDITION = {Second},
 PUBLISHER = {Springer-Verlag, Berlin},
      YEAR = {1998},
     PAGES = {xiv+470},
      ISBN = {3-540-62046-X; 0-387-98549-2},
   MRCLASS = {14C17 (14-02)},
  MRNUMBER = {1644323},
       DOI = {10.1007/978-1-4612-1700-8},
       URL = {https://doi.org/10.1007/978-1-4612-1700-8},
}

@book {Oxley,
    AUTHOR = {Oxley, James},
     TITLE = {Matroid theory},
    SERIES = {Oxford Graduate Texts in Mathematics},
    VOLUME = {21},
   EDITION = {Second},
 PUBLISHER = {Oxford University Press, Oxford},
      YEAR = {2011},
     PAGES = {xiv+684},
      ISBN = {978-0-19-960339-8},
   MRCLASS = {05-01 (05B35 90C27)},
  MRNUMBER = {2849819},
MRREVIEWER = {Maruti M. Shikare},
       DOI = {10.1093/acprof:oso/9780198566946.001.0001},
       URL = {https://doi.org/10.1093/acprof:oso/9780198566946.001.0001},
}

@article {JiangLi,
    AUTHOR = {Jiang, Chen and Li, Zhiyuan},
     TITLE = {Algebraic reverse {K}hovanskii-{T}eissier inequality via
              {O}kounkov bodies},
   JOURNAL = {Math. Z.},
  FJOURNAL = {Mathematische Zeitschrift},
    VOLUME = {305},
      YEAR = {2023},
    NUMBER = {2},
     PAGES = {Paper No. 26, 14},
      ISSN = {0025-5874,1432-1823},
   MRCLASS = {14C20 (14C17 14M25)},
  MRNUMBER = {4645768},
MRREVIEWER = {Pietro\ Sabatino},
       DOI = {10.1007/s00209-023-03349-9},
       URL = {https://doi.org/10.1007/s00209-023-03349-9},
}

@article {Yau,
    AUTHOR = {Yau, Shing Tung},
     TITLE = {On the {R}icci curvature of a compact {K}\"{a}hler manifold and
              the complex {M}onge-{A}mp\`ere equation. {I}},
   JOURNAL = {Comm. Pure Appl. Math.},
  FJOURNAL = {Communications on Pure and Applied Mathematics},
    VOLUME = {31},
      YEAR = {1978},
    NUMBER = {3},
     PAGES = {339--411},
      ISSN = {0010-3640},
   MRCLASS = {53C55 (32C10 35J60)},
  MRNUMBER = {480350},
MRREVIEWER = {Robert E. Greene},
       DOI = {10.1002/cpa.3160310304},
       URL = {https://doi.org/10.1002/cpa.3160310304},
}

@article {LehmannXiao,
    AUTHOR = {Lehmann, Brian and Xiao, Jian},
     TITLE = {Correspondences between convex geometry and complex geometry},
   JOURNAL = {\'{E}pijournal G\'{e}om. Alg\'{e}brique},
  FJOURNAL = {\'{E}pijournal de G\'{e}om\'{e}trie Alg\'{e}brique. EPIGA},
    VOLUME = {1},
      YEAR = {2017},
     PAGES = {Art. 6, 29},
   MRCLASS = {14C20 (32Q15 52A39)},
  MRNUMBER = {3743109},
MRREVIEWER = {Eugenii Shustin},
       DOI = {10.46298/epiga.2017.volume1.2038},
       URL = {https://doi.org/10.46298/epiga.2017.volume1.2038},
}

@article {Postnikov,
    AUTHOR = {Postnikov, Alexander},
     TITLE = {Permutohedra, associahedra, and beyond},
   JOURNAL = {Int. Math. Res. Not. IMRN},
  FJOURNAL = {International Mathematics Research Notices. IMRN},
      YEAR = {2009},
    NUMBER = {6},
     PAGES = {1026--1106},
      ISSN = {1073-7928},
   MRCLASS = {05E30},
  MRNUMBER = {2487491},
       DOI = {10.1093/imrn/rnn153},
       URL = {https://doi.org/10.1093/imrn/rnn153},
}

@book {FultonYoung,
    AUTHOR = {Fulton, William},
     TITLE = {Young tableaux},
    SERIES = {London Mathematical Society Student Texts},
    VOLUME = {35},
      NOTE = {With applications to representation theory and geometry},
 PUBLISHER = {Cambridge University Press, Cambridge},
      YEAR = {1997},
     PAGES = {x+260},
      ISBN = {0-521-56144-2; 0-521-56724-6},
   MRCLASS = {05E10 (05E05 05E15 14M15 20G05)},
  MRNUMBER = {1464693},
MRREVIEWER = {Tadeusz J\'{o}zefiak},
}

@incollection {GurvitsL,
    AUTHOR = {Gurvits, Leonid},
     TITLE = {On multivariate {N}ewton-like inequalities},
 BOOKTITLE = {Advances in combinatorial mathematics},
     PAGES = {61--78},
 PUBLISHER = {Springer, Berlin},
      YEAR = {2009},
   MRCLASS = {26D07 (05A20 26B25)},
  MRNUMBER = {2683227},
       DOI = {10.1007/978-3-642-03562-3\_4},
       URL = {https://doi.org/10.1007/978-3-642-03562-3_4},
}

@book {Murota,
    AUTHOR = {Murota, Kazuo},
     TITLE = {Discrete convex analysis},
    SERIES = {SIAM Monographs on Discrete Mathematics and Applications},
 PUBLISHER = {Society for Industrial and Applied Mathematics (SIAM),
              Philadelphia, PA},
      YEAR = {2003},
     PAGES = {xxii+389},
      ISBN = {0-89871-540-7},
   MRCLASS = {90-02 (52-02 90C27 90C46 91B02)},
  MRNUMBER = {1997998},
MRREVIEWER = {Ulrich Faigle},
       DOI = {10.1137/1.9780898718508},
       URL = {https://doi.org/10.1137/1.9780898718508},
}

@article {HMMS,
    AUTHOR = {Huh, June and Matherne, Jacob and M\'esz\'aros, Karola and
              St. Dizier, Avery},
     TITLE = {Logarithmic concavity of {S}chur and related polynomials},
   JOURNAL = {Trans. Amer. Math. Soc.},
  FJOURNAL = {Transactions of the American Mathematical Society},
    VOLUME = {375},
      YEAR = {2022},
    NUMBER = {6},
     PAGES = {4411--4427},
      ISSN = {0002-9947,1088-6850},
   MRCLASS = {05E10 (14M15 17B10)},
  MRNUMBER = {4419063},
MRREVIEWER = {Allan\ Berele},
       DOI = {10.1090/tran/8606},
       URL = {https://doi.org/10.1090/tran/8606},
}

@book {FultonToric,
    AUTHOR = {Fulton, William},
     TITLE = {Introduction to toric varieties},
    SERIES = {Annals of Mathematics Studies},
    VOLUME = {131},
 PUBLISHER = {Princeton University Press, Princeton, NJ},
      YEAR = {1993},
     PAGES = {xii+157},
      ISBN = {0-691-00049-2},
   MRCLASS = {14M25 (14-02 14J30)},
  MRNUMBER = {1234037},
MRREVIEWER = {T.\ Oda},
       DOI = {10.1515/9781400882526},
       URL = {https://doi.org/10.1515/9781400882526},
}

@article {Eilenberg,
    AUTHOR = {Eilenberg, Samuel},
     TITLE = {On the problems of topology},
   JOURNAL = {Ann. of Math. (2)},
  FJOURNAL = {Annals of Mathematics. Second Series},
    VOLUME = {50},
      YEAR = {1949},
     PAGES = {247--260},
      ISSN = {0003-486X},
   MRCLASS = {56.0X},
  MRNUMBER = {30189},
MRREVIEWER = {H.\ Cartan},
       DOI = {10.2307/1969448},
       URL = {https://doi.org/10.2307/1969448},
}

@book {Lazarsfeld1,
    AUTHOR = {Lazarsfeld, Robert},
     TITLE = {Positivity in algebraic geometry. {I}},
    SERIES = {Ergebnisse der Mathematik und ihrer Grenzgebiete},
    VOLUME = {48},
      NOTE = {Classical setting: line bundles and linear series},
 PUBLISHER = {Springer-Verlag, Berlin},
      YEAR = {2004},
     PAGES = {xviii+387},
      ISBN = {3-540-22533-1},
   MRCLASS = {14-02 (14C20)},
  MRNUMBER = {2095471},
MRREVIEWER = {Mihnea\ Popa},
       DOI = {10.1007/978-3-642-18808-4},
       URL = {https://doi.org/10.1007/978-3-642-18808-4},
}

@misc{BrandenLeake,
      title={Lorentzian polynomials on cones}, 
      author={Br\"and\'en, Petter and Leake, Jonathan},
      eprint={2304.13203},
      note={arXiv:2304.13203},
      url={https://arxiv.org/abs/2304.13203}, 
}

@article {Garding,
    AUTHOR = {Garding, Lars},
     TITLE = {Linear hyperbolic partial differential equations with constant
              coefficients},
   JOURNAL = {Acta Math.},
  FJOURNAL = {Acta Mathematica},
    VOLUME = {85},
      YEAR = {1951},
     PAGES = {1--62},
      ISSN = {0001-5962,1871-2509},
   MRCLASS = {36.0X},
  MRNUMBER = {41336},
MRREVIEWER = {F.\ John},
       DOI = {10.1007/BF02395740},
       URL = {https://doi.org/10.1007/BF02395740},
}

@article {Minkowski,
    AUTHOR = {Minkowski, Hermann},
     TITLE = {Volumen und {O}berfl\"ache},
   JOURNAL = {Math. Ann.},
  FJOURNAL = {Mathematische Annalen},
    VOLUME = {57},
      YEAR = {1903},
    NUMBER = {4},
     PAGES = {447--495},
      ISSN = {0025-5831,1432-1807},
   MRCLASS = {99-04},
  MRNUMBER = {1511220},
       DOI = {10.1007/BF01445180},
       URL = {https://doi.org/10.1007/BF01445180},
}

@Article{BrandenHuh,
 Author = {Br{\"a}nd{\'e}n, Petter and Huh, June},
 Title = {Lorentzian polynomials},
 FJournal = {Annals of Mathematics. Second Series},
 Journal = {Ann. Math. (2)},
 ISSN = {0003-486X},
 Volume = {192},
 Number = {3},
 Pages = {821--891},
 Year = {2020},
 DOI = {10.4007/annals.2020.192.3.4},
 zbMATH = {7285355},
 Zbl = {1454.52013},
  MRNUMBER = {4187216},
}

@article {RSW,
    AUTHOR = {Ross, Julius and S\"uss, Hendrik and Wannerer, Thomas},
     TITLE = {Dually {L}orentzian {P}olynomials},
   JOURNAL = {Monatsh. Math.},
  FJOURNAL = {Monatshefte f\"ur Mathematik},
    VOLUME = {208},
      YEAR = {2025},
    NUMBER = {3},
     PAGES = {495--524},
      ISSN = {0026-9255,1436-5081},
   MRCLASS = {14C17 (32J27 52A39 52A40 52B40)},
  MRNUMBER = {4994724},
       DOI = {10.1007/s00605-025-02134-6},
       URL = {https://doi.org/10.1007/s00605-025-02134-6},
}

@article {Aluffi,
    AUTHOR = {Aluffi, Paolo},
     TITLE = {Lorentzian polynomials, {S}egre classes, and adjoint
              polynomials of convex polyhedral cones},
   JOURNAL = {Adv. Math.},
  FJOURNAL = {Advances in Mathematics},
    VOLUME = {437},
      YEAR = {2024},
     PAGES = {Paper No. 109440, 37},
      ISSN = {0001-8708,1090-2082},
   MRCLASS = {14C17 (52A20)},
  MRNUMBER = {4674860},
MRREVIEWER = {I.\ Dolgachev},
       DOI = {10.1016/j.aim.2023.109440},
       URL = {https://doi.org/10.1016/j.aim.2023.109440},
}

@article {HuhChromatic,
    AUTHOR = {Huh, June},
     TITLE = {Milnor numbers of projective hypersurfaces and the chromatic
              polynomial of graphs},
   JOURNAL = {J. Amer. Math. Soc.},
  FJOURNAL = {Journal of the American Mathematical Society},
    VOLUME = {25},
      YEAR = {2012},
    NUMBER = {3},
     PAGES = {907--927},
      ISSN = {0894-0347,1088-6834},
   MRCLASS = {14B05 (05B35 14C17)},
  MRNUMBER = {2904577},
MRREVIEWER = {Paolo\ Aluffi},
       DOI = {10.1090/S0894-0347-2012-00731-0},
       URL = {https://doi.org/10.1090/S0894-0347-2012-00731-0},
}

@Article{Shephard,
 Author = {Shephard, Geoffrey C.},
 Title = {Inequalities between mixed volumes of convex sets},
 FJournal = {Mathematika},
 Journal = {Mathematika},
 ISSN = {0025-5793},
 Volume = {7},
 Pages = {125--138},
 Year = {1960},
 Language = {English},
 DOI = {10.1112/S0025579300001674},
 zbMATH = {3177259},
 Zbl = {0108.35203}
}

@InCollection{HuhICM,
 Author = {Huh, June},
 Title = {Combinatorics and {Hodge} theory},
 BookTitle = {International congress of mathematicians 2022, ICM 2022, Helsinki, Finland, virtual, July 6--14, 2022. Volume 1. Prize lectures},
 ISBN = {978-3-98547-059-4; 978-3-98547-559-9; 978-3-98547-058-7; 978-3-98547-558-2},
 Pages = {212--239},
 Year = {2023},
 Publisher = {Berlin: European Mathematical Society (EMS)},
 Language = {English},
 DOI = {10.4171/ICM2022/205},
 zbMATH = {7822675},
 Zbl = {1535.05058}
}

@misc{HHMWW,
author={Daoji Huang and June Huh and Mateusz Micha{\l}ek and Botong Wang and Shouda Wang},
title={Realizations of homology classes and projection areas},
    note={arXiv:2505.08881},
    eprint={2505.08881},
    archivePrefix={arXiv},
    primaryClass={math.AG}
}

@misc{Hochstattler,
author={Winfried Hochst\"attler},
title={A family of rank $4$ non-algebraic matroids with pseudomodular dual},
    note={arXiv:2511.10417},
    eprint={2511.10417},
    archivePrefix={arXiv},
    primaryClass={math.CO}
}

@article {BFP,
    AUTHOR = {Bamiloshin, Michael and Farr\`as, Oriol and Padr\'o, Carles},
     TITLE = {A note on extension properties and representations of
              matroids},
   JOURNAL = {Discrete Appl. Math.},
  FJOURNAL = {Discrete Applied Mathematics. The Journal of Combinatorial
              Algorithms, Informatics and Computational Sciences},
    VOLUME = {376},
      YEAR = {2025},
     PAGES = {270--280},
      ISSN = {0166-218X,1872-6771},
   MRCLASS = {05B35},
  MRNUMBER = {4926337},
MRREVIEWER = {James\ Dylan\ Douthitt},
       DOI = {10.1016/j.dam.2025.06.028},
       URL = {https://doi.org/10.1016/j.dam.2025.06.028},
}

@misc{BHMPW,
author={Tom Braden and June Huh and Jacob Matherne and Nicholas Proudfoot and Botong Wang},
title={Singular {H}odge theory for combinatorial geometries},
    note={arXiv:2010.06088},
    eprint={2010.06088},
    archivePrefix={arXiv},
    primaryClass={math.AG}
}

@misc{HMWX,
author={June Huh and Mateusz Micha{\l}ek and Botong Wang and Jian Xiao},
title={Correlation inequalities for volume polynomials},
    note={in preparation},
    eprint={},
    archivePrefix={arXiv},
    primaryClass={math.AG}
}

@misc{GHMSW,
author={Lukas Grund and June Huh and Mateusz Micha{\l}ek and Hendrik S{\"u}{\ss} and Botong Wang},
title={Linear operators preserving volume polynomials},
    note={arXiv:2506.22415},
    eprint={2506.22415},
    archivePrefix={arXiv},
    primaryClass={math.AG}
}

@article{Heine,
  author       = {Rudolf Heine},
  title        = {Der Wertvorrat der gemischten Inhalte von zwei, drei und vier ebenen Eibereichen},
  journal      = {Math. Ann.},
  volume       = {115},
  year         = {1938},
  number       = {1},
  pages        = {115--129},
  mrnumber     = {1513176},
  language     = {German},
}

@misc{HuhCorrespondence,
author={June Huh},
title={Correspondences between projective planes},
    note={arXiv:1303.4113},
    eprint={1303.4113},
    archivePrefix={arXiv},
    primaryClass={math.AG}
}

@book {WelshMatroidTheory,
    AUTHOR = {Welsh, D. J. A.},
     TITLE = {Matroid theory},
    SERIES = {L. M. S. Monographs},
    VOLUME = {No. 8},
 PUBLISHER = {Academic Press [Harcourt Brace Jovanovich, Publishers],
              London-New York},
      YEAR = {1976},
     PAGES = {xi+433},
   MRCLASS = {05B35},
  MRNUMBER = {427112},
MRREVIEWER = {W.\ T.\ Tutte},
}

@article {Lindstrom89,
    AUTHOR = {Lindstr\"om, B.},
     TITLE = {Matroids algebraic over {$F(t)$} are algebraic over {$F$}},
   JOURNAL = {Combinatorica},
  FJOURNAL = {Combinatorica. An International Journal of the J\'anos Bolyai
              Mathematical Society},
    VOLUME = {9},
      YEAR = {1989},
    NUMBER = {1},
     PAGES = {107--109},
      ISSN = {0209-9683},
   MRCLASS = {05B35 (12F20)},
  MRNUMBER = {1010307},
MRREVIEWER = {Dragan\ M.\ Acketa},
       DOI = {10.1007/BF02122691},
       URL = {https://doi.org/10.1007/BF02122691},
}

@misc{BHKL-triangular,
      title={Lorentzian polynomials and matroids over triangular hyperfields 1: Topological aspects}, 
      author={Matthew Baker and June Huh and Mario Kummer and Oliver Lorscheid},
      note= {arXiv:2508.02907},
      eprint={arXiv:2508.02907},
      archivePrefix={arXiv},
      primaryClass={math.CO},
      url={https://arxiv.org/abs/2508.02907}, 
}

@misc{BHKL-triangular2,
      title={Lorentzian polynomials and matroids over triangular hyperfields 2: Analytic aspects}, 
      author={Matthew Baker and June Huh and Mario Kummer and Oliver Lorscheid},
      note={in preparation},
      eprint={},
      archivePrefix={},
      primaryClass={math.CO},
      url={}, 
}

@misc{BHKL-polymatroids,
      title={Representation theory for polymatroids}, 
      author={Matthew Baker and June Huh and Donggyu Kim and Mario Kummer and Oliver Lorscheid},
      note= {arXiv:2507.14718},
      eprint={2507.14718},
      archivePrefix={arXiv},
      primaryClass={math.CO},
      url={https://arxiv.org/abs/2507.14718}, 
}

\end{document}